\documentclass[12pt]{amsart}

\usepackage{latexsym}

\newtheorem{thm}{Theorem}[section]
\newtheorem{defn}[thm]{Definition}
\newtheorem{lem}[thm]{Lemma}%[section]
\newtheorem{prop}[thm]{Proposition}%[section]
\newtheorem{cor}[thm]{Corollary}%[section]

\newenvironment{rem}{\medskip\noindent{\it Remark:\/} }{\medskip}

\newcommand{\szego}{Szeg\"o }

\newcommand{\kahler}{K\"ahler }

\newcommand{\R}{{\mathbb R}}
\newcommand{\C}{{\mathbb C}}

\newcommand{\Z}{{\mathbb Z}}

\newcommand{\dbar}{\bar\partial}
\newcommand{\ddbar}{\partial\dbar}

\renewcommand{\H}{{\mathbf H}}

\newcommand{\half}{{\frac{1}{2}}}

\renewcommand{\phi}{\varphi}

\newcommand{\hcal}{\mathcal{H}}

\newcommand{\lcal}{\mathcal{L}}

\newcommand{\al}{\alpha}

\newcommand{\la}{\lambda}

\begin{document}

\title
{Quantum maps and automorphisms}
\author{Steve Zelditch }
\address{Department of Mathematics, Johns Hopkins University, Baltimore, MD
21218, USA} \email{zelditch@math.jhu.edu}

\thanks{Research partially supported by NSF grant  DMS-0071358 and by the Clay Mathematics Institute .}

\begin{abstract} There are several inequivalent definitions of what it means to
quantize a symplectic map on a symplectic manifold $(M, \omega)$.
One definition is that the quantization is an automorphism of a
$*$ algebra associated to $(M, \omega)$. Another is that it is
unitary operator $U_{\chi}$ on a Hilbert space associated to $(M,
g)$, such that $A \to U_{\chi}^* A U_{\chi}$ defines an
automorphism of the algebra of observables. A yet stronger one,
common in partial differential equations, is that $U_{\chi}$
should be a Fourier integral operator associated to the graph of
$\chi$. We compare the definitions in the case where $(M,\omega)$
is a compact \kahler manifold. The main result is a Toeplitz
analogue of the Duistermaat-Singer theorem on automorphisms of
pseudodifferential algebras, and an extension which does not
assume $H^1(M, \C) = \{0\}$.  We illustrate with examples from
quantum maps.
\end{abstract}

\maketitle

\section{Introduction}

Much attention has been focussed recently on  $*$ products on
Poisson manifolds $(M, \{, \})$ (see, among others, \cite{ K1, CF,
KS,  T, RT, EK, WX} ). Such $*$ products are  viewed as quantizing
functions on $M$ to an algebra of observables.  This article is
concerned with the related problem  of quantizing symplectic maps
 $\chi$ on K\"ahler manifolds $(M, \omega)$,  a special case of the problem of quantizing Poisson maps.
  From the * algebra viewpoint,  it seems most natural to quantize such a symplectic map as an automorphism of a  $*$ algebra
  associated to $(M, \omega)$, specifically
 the (complete)
 symbol algebra   ${\mathcal
T}^*/{\mathcal T}^{-\infty}$  of Berezin-Toeplitz operators
 over   $(M, \omega)$.
These symbol  algebras  are  basic examples of  abstract $*$
algebras arising in deformation quantization of  Poisson manifolds
(see \cite{Bou1, Bou2, Bou3, Ch, CF, Gu, Sch1, Sch2} for more on
this aspect). But they carry more structure than bare  * algebras:
the Toeplitz operator algebra ${\mathcal T}^*$ of which ${\mathcal
T}^*/{\mathcal T}^{-\infty}$ is the symbol algebra also comes with
a representation as operators on a Hilbert space. In the Hilbert
space setting, it  is most natural to try to quantize a symplectic
map $\chi$ as a unitary operator $U_{\chi}$, and to induce the
automorphism $U_{\chi} A U_{\chi}^*$ on ${\mathcal T}^*$ . As will
be explained below (see also \cite{Zel1}), it is not always
possible to quantize a symplectic map this way. When possible, the
quantization $U_{\chi}$ is an example of what is known as a
quantum map   in the literature of quantum chaos. Such quantum
maps have also been the focus of much attention in recent years by
a virtually disjoint group  (see e.g. \cite{dBdE, Kea, deGI, HB,
MR, Zel1}). The main purpose of this article is to contrast the
different notions of quantizing symplectic maps, as  as they arise
in Toeplitz * algebras, partial differential equations and quantum
chaos.

 Aside
from its intrinsic interest, the relation between quantum maps and
automorphisms of * algebras  has practical consequences in quantum
chaos, i.e. in the relations between dynamical properties of
$\chi$ and the eigenvalues/eigenfunctions of its quantization
$U_{\chi}$. In the physics literature of quantum chaos, quantum
maps are  studied through
 examples such  as   quantum kicked  tops (on $S^2$),   cat maps,
 rotors, baker's map and standard maps (on the 2- torus ${\bf T}^2$).
Almost always, the quantizations are given  as explicit unitary
matrices $U_{ N}$  (depending on a Planck constant $1/N$) on
special Hilbert spaces ${\mathcal H}_N$, often using some special
representation theory, and no formal definition is given of the
term `quantum map'. The need for precise definitions  is felt,
however, as soon as one aims at quantizing maps which lie outside
the range of standard examples. Even symplectic maps on surfaces
of genus $g \geq 2$ count as non-standard, and only seem to have
been quantized by the Toeplitz method discussed in this paper and
in  \cite{Zel1}.

A further reason to study quantum maps versus automorphisms is to
better understand obstructions to quantizations. It is often said
that Kronecker translations
$$ T_{\alpha, \beta} (x, \xi) =  (x + \alpha, \xi + \beta), \;\;\; (x, \xi) \in \R^{2n}/\Z^{2n},$$
and  affine symplectic torus  maps
$$ f_{\alpha}(x, \xi) = (x + \xi, \xi + \alpha)$$
are  not quantizable, for  reasons explained in Proposition
\ref{SYMPLIFT}. Nevertheless, the paper  \cite{MR}
 proposes a quantization of such maps. Of course, the resolution
 of this paradox is that a weaker notion of quantization is
 assumed in \cite{MR} than elsewhere, as will be explained below.

The implicit criterion (including that in  \cite{MR})  that $U_N$
quantize a symplectic map $\chi$ is that the Egorov type formula
\begin{equation} \label{EGOROV} U_{ N}^* Op_N(a) U_{ N} \sim  Op_N(a \circ
\chi),\;\; (N \to \infty) \end{equation} hold for all elements
$Op_N(a) $ of the algebra ${\mathcal T}_N$ of observables, where
$\chi$ is a symplectic map of $(M, \omega)$. Postponing precise
definitions, we see that the operative condition is that $U_{ N}^*
Op_N(a) U_{ N}$ defines an automorphism  of ${\mathcal T}_N$, at
least to leading order. Here, our notation for  observables and
quantum maps are in terms of  sequences as the inverse Planck
constnat  $N$ varies. We temporarily write ${\mathcal T}^*$ for
sequences $\{Op_N(a)\}$ of observables (with ${\mathcal
T}^{-\infty}$ the sequences which are rapidly decaying in $N$),
and $U \sim \{U_N\}$ for sequences of unitary quantum maps. We
will soon give more precise definitions.

We now  distinguish several notions of quantizing a symplectic map
and make a number of assertions which will be justified in the
remainder of the article.
\begin{itemize}

\item   There is a geometric obstruction to quantizing a symplectic map
$\chi$ as a Toeplitz quantum map $U_{\chi, N}$ on ${\mathcal H}_N$
(see Definition \ref{TQM} and Proposition \ref{SYMPLIFT}) .
Kronecker translations and parabolic maps of the torus are
examples of non-quantizable symplectic maps in the Toeplitz sense
(see Propositions \ref{KRON} and \ref{PARA});

\item  There is no obstruction to quantizing a symplectic map as an
automorphism of the Toeplitz symbol algebra  ${\mathcal
T}^*/{\mathcal T}^{-\infty}$ (see Theorem \ref{DSMEETTOP}). For
instance, Kronecker maps and cat maps are quantizable as
automorphisms (see Propositions \ref{KRON2}-\ref{PARA2});

\item  Conversely, if $H^1(M, \C) = \{0\}$, then every order preserving  automorphism of the
symbol algebra ${\mathcal T}^*/{\mathcal T}^{-\infty}$ on $M$ is
induced by a symplectic map of $(M, \omega)$ (see Theorem
\ref{DSMEETTOP} for this and for the case where $H^1(M, \C) \not=
\{0\} );$

\item  There  is an obstruction to `extending' an automorphism
$\alpha$
 of ${\mathcal T}^*/{\mathcal T}^{-\infty}$
as an automorphism of ${\mathcal T}^*$. In particular, there is an
obstruction to inducing
  automorphisms $\alpha_N$ of the finite dimensional
algebras of operators ${\mathcal T}_N$ acting on ${\mathcal H}_N$
(see Theorem \ref{DSMEETTOP}). Again, Kronecker maps are examples
(see \S \ref{QTM}).

\item Any sequence $U_N$ of unitaries on ${\mathcal H}_N$ which
defines an automorphism of ${\mathcal T}^*/{\mathcal T}^{-\infty}$
must be a Toeplitz quantum map in sense of Definition \ref{TQM}
(i) (cf. \cite{Bou4, Zel1}).

\item  Many of the key problems of quantum chaos, e.g. problems on
eigenvalue  level spacings or pair correlation, on ergodicity and
mixing of eigenfunctions (etc.)  concern
 only the spectral theory of the  automorphism quantizing $\chi$ and
not the unitary map per se (see \S \ref{SP}).

\end{itemize}

\subsection{The Toeplitz set-up}

In order to state our results precisely, we need to specify the
framework in which we are working. The framework of Toeplitz
operators  used in this paper is the same as in \cite{Bou1, Bou2,
Gu,  BSZ, SZ, Zel1, Zel2}. We briefly recall the notation and
terminology.

Our setting consists of a  \kahler manifold $(M, \omega)$ with
 $\frac{1}{2\pi} [\omega] \in H^1(M, \Z)$. Under this integrality condition,
 there exists a     positive
 hermitian holomorphic
line bundle $(L, h) \to M$ over $M$ with  curvature form
$$c_1(h)=-\frac{\sqrt{-1}}{\pi}\ddbar \log \|e_L\|_h\; = \omega,$$ where $e_L$ is a
nonvanishing local holomorphic section of $L$, and where
$\|e_L\|_h=h(e_L,e_L)^{1/2}$ denotes the $h$-norm of $e_L$. We
give $M$ the volume form $ dV= \frac{1}{m!} \omega^m.$

The Hilbert spaces `quantizing' $(M, \omega)$  are then defined to
be the  spaces $H^0(M,L^N)$ of holomorphic sections of
$L^N=L\otimes\cdots\otimes L$.  The metric $h$ induces Hermitian
metrics $h_N$ on $L^N$ given by $\|s^{\otimes
N}\|_{h_N}=\|s\|_h^N$.  We give $H^0(M,L^N)$ the inner product
\begin{equation}\label{inner}\langle s_1, s_2 \rangle = \int_M h_N(s_1,
s_2)dV \quad\quad (s_1, s_2 \in H^0(M,L^N)\,)\;,\end{equation} and
we write $|s|=\langle s,s \rangle^{1/2}$. We then define the
Szeg\"o kernels as the orthogonal projections $\Pi_N: \lcal^2(M,
L^N) \to H^0(M, L^N)$, so that
\begin{equation} (\Pi_N s)(w)  = \int_M h^N_z\big(s(z),\Pi_N(z,w)\big)
dV_M(z)\,, \quad s\in \lcal^2(M, L^N)\,.
\end{equation}

Instead of dealing with sequences of Hilbert spaces, observables
and unitary operators, it is convenient to lift them to  the
circle bundle $X=\{\la \in L^* : \|\la\|_{h^*}= 1\}$, where $L^*$
is the dual line bundle to $L$, and  where $h^*$ is the norm on
$L^*$ dual to $h$. Associated to $X$ is the contact form $\al=
-i\partial\rho|_X=i\dbar\rho|_X$ and  the volume form
\begin{equation}\label{dvx}dV_X=\frac{1}{m!}\al\wedge
(d\al)^m=\al\wedge\pi^*dV_M\,.\end{equation}

Holomorphic sections then lift to elements of the  Hardy space
$\hcal^2(X) \subset \lcal^2(X)$ of square-integrable CR functions
on $X$, i.e., functions that are annihilated by the Cauchy-Riemann
operator $\dbar_b$ and are $\lcal^2$ with respect to the inner
product
\begin{equation}\label{unitary} \langle  F_1, F_2\rangle
=\frac{1}{2\pi}\int_X F_1\overline{F_2}dV_X\,,\quad
F_1,F_2\in\lcal^2(X)\,.\end{equation}  We let $r_{\theta}x
=e^{i\theta} x$ ($x\in X$) denote the $S^1$ action on $X$ and
denote its infinitesimal generator by $\frac{\partial}{\partial
\theta}$. The $S^1$ action on $X$ commutes with
$\bar{\partial}_b$; hence $\hcal^2(X) = \bigoplus_{N =0}^{\infty}
\hcal^2_N(X)$ where $\hcal^2_N(X) = \{ F \in \hcal^2(X):
F(r_{\theta}x) = e^{i N \theta} F(x) \}$. A section $s_N$ of $L^N$
determines an equivariant function $\hat{s}_N$ on $L^*$ by the
rule
$$\hat{s}_N(\lambda) = \left( \lambda^{\otimes N}, s_N(z)
\right)\,,\quad \la\in L^*_z\,,\ z\in M\,,$$ where
$\lambda^{\otimes N} = \lambda \otimes \cdots\otimes \lambda$. We
henceforth restrict $\hat{s}$ to $X$ and then the equivariance
property takes the form $\hat s_N(r_\theta x) = e^{iN\theta} \hat
s_N(x)$. The map $s\mapsto \hat{s}$ is a unitary equivalence
between $H^0(M, L^{ N})$ and $\hcal^2_N(X)$. We refer to \cite{BG,
BS,  BSZ, Zel2} for further background.

We now define the (lifted) \szego kernel of degree $N$ to be the
orthogonal projection $\Pi_N : \lcal^2(X)\rightarrow
\hcal^2_N(X)$. It is defined by
\begin{equation} \Pi_N F(x) = \int_X \Pi_N(x,y) F(y) dV_X (y)\,,
\quad F\in\lcal^2(X)\,. \label{PiNF}\end{equation} The full \szego
kernel is the direct sum
\begin{equation} \Pi = \bigoplus_{N=1}^{\infty} \Pi_N. \end{equation}
 Following Boutet de Monvel-Guillemin \cite{BG}, we then define:

 \begin{defn} The  * algebra  ${\mathcal T}^*(M)$ of Toeplitz operators of $(M, \omega)$ is the algebra
 of operators on $\hcal^2(X)$ of the form \begin{equation}  \Pi A \Pi =
\bigoplus_{N=1}^{\infty} \Pi_N A_N \Pi_N,\;\; A \in \Psi^*_{S^1}
(X)
\end{equation} where $\Psi^*_{S^1}(X)$ is the algebra of pseudodifferential operators over $X$
which commute with the $S^1$ action, and where
\begin{equation} A_N = \int_{S^1} e^{i {\mathcal N} \theta} A e^{- i {\mathcal N} \theta}  d \theta. \end{equation}
\end{defn}
Here,  ${\mathcal N} = \frac{1}{i} \frac{\partial }{\partial
\theta}$ is the operator generating the $S^1$ action, whose
eigenvalue in $\hcal^2_N(X)$ equals $N$.

Since the symbol of $A$ is $S^1$-invariant,  Toeplitz operators of
this kind possess an expansion
\begin{equation} \label{SPECIAL} \Pi A \Pi \sim {\mathcal N}^s \;
\sum_{j = 0}^{\infty}\Pi  a_j \Pi {\mathcal N}^{-j} \end{equation}
where $a_j \in C^{\infty}(M)$. We may also express it
 in the direct sum form
\begin{equation} \Pi A \Pi = \sum_{N = 1}^{\infty} \Pi_N a_N \Pi_N \end{equation}
where $a_N(z, \bar{z}) \in S_{scl}^{s} $ is a semiclassical symbol
of some order $s$, i.e. admits an asymptotic expansion
\begin{equation} \label{SC} a_N(z, \bar{z}) \sim N^{s} \sum_{j = 0}^{\infty} N^{-j} a_j(z, \bar{z}),\;\;\;\;
a_j(z, \bar{z}) \in C^{\infty}(M)
\end{equation}
in the sense of symbols. We define the {\it order} of a Toeplitz
operator $\Pi A \Pi$ to be the order $s$  of the  symbol. The
order defines a filtration of ${\mathcal T}^*$  by spaces of
operators ${\mathcal T}^s$ of order $s \in \R.$ See \cite{Gu} for
further background.

We also define `flat' symbols $f(z,N) \in S_{scl}^{- \infty}$ is
$\sim 0$ as functions satisfying $f = O(N^{-m})$ for all $m$. We
then define  ${\mathcal T}^{-\infty}$ to be the flat (or
smoothing) Toeplitz operators (possessing a flat symbol).  The
following definition is important in distinguishing the
automorphisms which concern us:
\begin{defn} The  complete Toeplitz symbol algebra (or smooth  Toeplitz algebra)  is the quotient algebra
 $\mathcal
T^{*}/\mathcal T^{- \infty}$.
 \end{defn}

We often view $\bigoplus_{N = 1}^{\infty} \Pi_N a_N \Pi_N$ as the
sequence $\{\Pi_N a_N \Pi_N\}$ of operators on the sequence
$\hcal_N \simeq H^0(M, L^N)$ of Hilbert spaces.  The physicists'
notation for $\Pi_N a_N \Pi_N$ is $Op_N(a_N)$.  Viewing symbols as
sequences $\{a_N(z,\bar{z})$, we define the $*_N$ product by
\begin{equation}\label{STAR}  \Pi_N a_N \Pi_N \circ \Pi_N b_N \Pi_N = \Pi_N a_N
*_N b_N \Pi_N.
\end{equation} In the Appendix, we will describe the calculation
of $a_N *_N b_N$ so that it will not seem abstract to the reader.
We now introduce automorphisms:

\begin{defn}  An  order preserving automorphism $\alpha$ of ${\mathcal T}^*/{\mathcal T}^{-\infty}$ is an
automorphism  which preserves the filtration ${\mathcal
T}^s/{\mathcal T}^{-\infty}$. We denote the algebra of such
automorphisms by $Aut_o({\mathcal T}^*/{\mathcal T}^{-\infty})$.

\end{defn}

It is important to distinguish:
\begin{itemize}

\item Order preserving automorphisms of ${\mathcal T}^*$ which preserve ${\mathcal
T}^{-\infty}$;

\item Order preserving automorphisms of the symbol algebra ${\mathcal T}^*/{\mathcal
T}^{-\infty}$.

\end{itemize}

Since elements of ${\mathcal T}^*$ and of  ${\mathcal
T}^*/{\mathcal T}^{-\infty}$ commute with the $S^1$ action, either
kind of automorphism satisfies:
\begin{equation} \alpha(\bigoplus_{N = 1}^{\infty} \Pi_N a_N
\Pi_N) \sim  \bigoplus_{N = 1}^{\infty} \Pi_N b_N \Pi_N,
\end{equation}
where $b_N$ is a semiclassical symbol of the same order as $a_N$.
In the case of automorphisms of ${\mathcal T}^*$,  we can conclude
that $\alpha (\Pi_N a_N \Pi_N) = \Pi_N b_N \Pi_N$ and that
$\alpha$ induces automorphisms $\alpha_N$ of the finite
dimensional algebras ${\mathcal T}_N$ for fixed $N$. However, for
automorphisms of ${\mathcal T}^*/{\mathcal T}^{-\infty}$ in
general, $\alpha (\Pi_N a_N \Pi_N)$ is not even defined since
 $\Pi_N a_N \Pi_N \in
{\mathcal T}^{-\infty}$. To put it another way, we cannot uniquely
represent an element of the finite dimensional algebra as $\Pi_N
a_N \Pi_N$ although we can uniquely represent elements of
${\mathcal T}^*/{\mathcal T}^{-\infty}$ this way.

\subsubsection{\label{COCO} Covariant and Contravariant symbols}

Let $\Pi_N a \Pi_N$ be a Toeplitz operator. By the {\it
contravariant symbol} of $\Pi_N a_N \Pi_N$ is meant the multiplier
$a_N$.  By the {\it covariant symbol} of an operator $F$ is meant
the function \begin{equation} \hat{f}(z, \bar{z}) = \frac{\langle
F \Phi_N^z, \Phi_N^w \rangle} { \langle  \Phi_N^z, \Phi_N^w
\rangle} |_{z = w} = \frac{\Pi_N F \Pi_N(z,z)}{\Pi_N(z,z)}.
\end{equation} where \begin{equation} \Phi_N^w(z) =
\frac{\Pi_N(z,w)}{||\Pi_N(\cdot, z)||} \end{equation}
 is
 the $L^2$-normalized  `coherent state' centered at $w$.
When $F = \Pi_N a \Pi_N$ we get
\begin{equation} \hat{a}(z, \bar{z}) = \frac{\Pi_N a \Pi_N(z,z)}{\Pi_N(z,z)}. \end{equation}
We use the notation $ I_N(a) = \hat{a}$ for the linear operator
(the Berezin transform)  which takes the contravariant symbol to
the covariant symbol (see \cite{RT} for background).

\subsection{Statement of results}

 Let us now consider the
senses in which we can quantize symplectic maps in our setting.
The first sense is that of  quantizations of symplectic maps as
Toeplitz Fourier integral operators. The definition is as follows.

\begin{defn} \label{TFIO} Suppose that the symplectic map $\chi$
of $(M, \omega)$ lifts to $(X, \alpha) $ as a contact
transformation $\tilde{\chi}$. By the Toeplitz Fourier integral
operator (or quantum map) defined by $\chi$ we mean the operator,
$$U = \bigoplus_N U_{\chi, N},\;\; U_{\chi, N} = \Pi_N T_{\chi} \sigma_N \Pi_N$$ where $T_{\chi}:
\lcal^2(X) \to \lcal^2(X)$ is the translation $T_{\chi}(f) = f
\circ \tilde{\chi}^{-1}$ and where $\sigma_N$ is a symbol designed
to make $U_{\chi, N}$ unitary. (Such a symbol always  exists
\cite{Zel1}).
\end{defn}

We now distinguish several notions of quantizing a symplectic map.

\begin{defn} \label{TQM} Let $\chi$ be a symplectic map of $(M, \omega)$. In
descending strength, we say that:

\begin{itemize}

\item  (a) $\chi$ is quantizable as a Toeplitz quantum map (or Toeplitz Fourier integral operator)  if it lifts
to a contact transformation $\tilde{\chi}$ of $(X, \alpha).$ The
quantization is then that of Definition \ref{TFIO};

\item  (b) $\chi$ is quantizable as an automorphism of the
full observable algebra if there exists an automorphism $\alpha$
of ${\mathcal T}^*$ satisfying (\ref{EGOROV});

\item (c)  $\chi$ is quantizable as an automorphism of the
symbol  algebra if there exists an automorphism $\alpha$  of
${\mathcal T}^*/{\mathcal T}^{-\infty}$ satisfying (\ref{EGOROV});

\end{itemize}

\end{defn}

By descending strength, we mean that quantization in a sense above
implies quantization in all of the following senses. The
automorphisms above are {\it order-preserving} in the sense that
the order of $\alpha(\Pi A \Pi)$ is the same as the order of $\Pi
A \Pi$. Henceforth, all automorphisms will be assumed to be
order-preserving.

We now explain the relations between these notions of
quantization.  We are guided in part by the analogous relations
between quantizations of symplectic maps (of cotangent bundles)
and automorphisms of the symbol algebra $\Psi^*/\Psi^{-\infty}$ of
the  algebra of pseudodifferential operators, as  determined   by
Duistermaat-Singer in \cite{DS1, DS2}. Their main result was that,
if $H^1(S^*M, \C) = \{0\}$, then every order preserving
automorphism of $\Psi^*/\Psi^{-\infty}$ is either conjugation by
an elliptic Fourier integral operator associated to the symplectic
map or a transmission. We prove an analogous theorem for Toeplitz
operators and  also extend it to the case where the phase space is
not simply connected.

To state the results, we need some notation. We denote the
universal cover of  $(M, \omega)$ by  $\tilde{M}$ and denote the
group of deck transformations of the natural cover $p: \tilde{M}
\to M$ by $\Gamma$. We lift all objects on $M$ to $\tilde{M}$
under $p$. We denote by $T_{\gamma}$ the unitary operator of
translation by $\gamma$ on $\lcal^2(\tilde{M})$. We also denote by
${\mathcal T}^*_{\Gamma}$ the algebra of  $\Gamma$-invariant
Toeplitz operators on $\tilde{M}$. It is important to understand
that ${\mathcal T}^*_{\Gamma}$ is not isomorphic to the algebra of
Toeplitz operators on $M$ since there are non-trivial (smoothing)
operators which act trivially on automorphic (periodic) functions.
In other words, the representation of ${\mathcal T}^*_{\Gamma}$ on
automorphic sections has a kernel, which we denote by ${\mathcal
K}_{\Gamma}$, and  ${\mathcal T}^*(X) \simeq {\mathcal
T}^*_{\Gamma}(\tilde{X})/{\mathcal K}_{\Gamma}$. Automorphisms
which descend to the finite Toeplitz algebras are precisely those
which preserve the subalgebra ${\mathcal K}_{\Gamma}$. For further
discussion,  we refer to \S \ref{NOTSC}.

\begin{thm}\label{DSMEETTOP} With the above notation, we have:
\begin{itemize}

\item (0) (Essentially known) A symplectic map of $(M, \omega)$
lifts to a contact transformation of $(X, \alpha)$ and hence
defines a Toeplitz quantum map if and only if it preserves
holonomies of all closed curves of $M$.  (See Proposition
\ref{SYMPLIFT} of  \S \ref{HOL}).

\item (i) Any symplectic map of any compact \kahler manifold $(M, \omega)$,
is quantizable as an  automorphism of the algebra  ${\mathcal
T}^*/{\mathcal T}^{-\infty}$ of smooth Toeplitz operators over
$M$;

\item (ii) Suppose that $H^1(M, \C) = \emptyset.$ Then any order-preserving automorphism of ${\mathcal T}^*/{\mathcal T}^{-\infty}$
is given by conjugation with a Toeplitz Fourier integral operator
on $M$  associated to a symplectic map $\chi$ of $(M, \alpha)$.
(The map  lifts to a  contact transformation of $(X, \alpha)$ by
(0)).

\item (iii) Suppose $H^1(M, \C) \not= \emptyset.$ Then to each  automorphism of  ${\mathcal T}^*/{\mathcal T}^{-\infty}$
there corresponds a symplectic map $\chi$ of $(M, \omega)$  and a
Toeplitz Fourier integral operator  (Definition \ref{TQM})
$U_{\chi}$ on the universal cover $\tilde{M}$ which satisfy
$T_{\gamma}^* U_{\chi} T_{\gamma} = M_{\gamma} U_{\chi}$, where
$M_{\gamma}$ is a central operator. The automorphism $A \to
U_{\chi}^* A U_{\chi}$ is $\Gamma$-invariant, and defines an
order-preserving automorphism of the algebra ${\mathcal
T}^*_{\Gamma}$ which  induces $\alpha$ on the $\Gamma$-invariant
symbol algebra  ${\mathcal T}^{*}/{\mathcal T}^{- \infty}$.

\item  (iv) Let ${\mathcal K}_{\Gamma} = \ker \rho_{\Gamma},$ where
$\rho_{\Gamma}$ is the representation of ${\mathcal T}^*_{\Gamma}$
on $\Gamma$-automorphic functions on $\tilde{M}.$  If $\alpha$
preserves ${\mathcal K}_{\Gamma}$, then it induces an
order-preserving automorphism on ${\mathcal T}^*(M)$ and hence on
the finite rank observables $Op_N(a)$ on ${\mathcal H}_N$.

\end{itemize}

\end{thm}

We separate the proof into the cases $H^1(M, \C) = \{0\}$ in \S
\ref{SC}  and $H^1(M, \C) \not= \{0\}$ in \S \ref{NOTSC}. The
latter case is very common in the physics literature on quantum
maps. The difference between order preserving automorphisms of
${\mathcal T}^*$ and ${\mathcal T}^*/{\mathcal T}^{-\infty}$ is
very significant, and only the former automorphisms are quantum
maps in the physics sense. For instance, as  will be seen in $\S
2$, Kronecker maps and affine symplectic maps are quantizable as
automorphisms of the symbol algebra, but do not lift to contact
transformations of $X$, do not preserve the kernel ${\mathcal
K}_{\Gamma}$ and are therefore not automorphisms of ${\mathcal
T}^*.$ Regarding (iv), it is not clear to us whether this operator
condition is equivalent to the holonomy-preservation condition in
(0).

As a corollary, we prove a result which indicates that  the
physicists' quantum maps are necessarily Toeplitz quantum maps
once they are conjugated to the complex (Bargmann) picture.

\begin{cor} \label{PHYSICS} If  $\{U_N\}$ is a sequence of unitary operators on
${\mathcal H}_N$ and if $\bigoplus U_N$ defines an order
preserving automorphism $\alpha$ of ${\mathcal T}^*/{\mathcal
T}^{-\infty}$, then $\{U_{N}\}$ must be a Toeplitz Fourier
integral operator associated to a quantizable symplectic map in
the sense of Definition \ref{TFIO} .
\end{cor}

As a gauge of our definitions, let us reconsider the
Marklof-Rudnick quantizations mentioned above of Kronecker maps,
parabolic maps and other `non-quantizable' maps \cite{MR}. They
define a sequence $\{U_N\}$ of unitary operators on ${\mathcal
H}_N$ satisfying  the leading order condition (\ref{EGOROV}), but
not to any lower order.  Hence,  conjugation $U_{\chi} Op_N(a)
U_{\chi}^*$ of an observable in ${\mathcal T}^*$ by  their quantum
map is no longer an observable, i.e. it is not an element of
${\mathcal T}^*$. Rather, its  Toeplitz symbol  only possesses a
one term asymptotic expansion and is not a classical symbol. Hence
it need not correspond to a quantizable symplectic map.

In addition to Theorem \ref{DSMEETTOP}, we discuss a related issue
revolving around the quantum maps versus automorphisms
distinction:  From the viewpoint of quantum chaos, the main
interest in the quantum maps $U_{\chi, N}$ lies in their spectral
theory and its relation to the dynamics of $\chi$. This is only
well-defined when the associated symplectic map is quantizable in
the strong sense as a sequence of unitary operators on $\hcal_N$.
As stated in the corollary, the symplectic map must then lift to a
contact transformation. In the last section \S \ref{SP}, we point
out that even when the symplectic map is quantizable as a unitary
operator, it is often the automorphism it induces which is most
significant in quantum chaos. That is,  much of the spectral
theory in quantum chaos concerns the spectrum of the automorphism
induced by $U_{\chi, N}$ rather than the spectrum of $U_{\chi, N}$
itself.

 The author benefited from discussions with S. de
Bievre,  Z. Rudnick  and particularly S. Nonnenmacher, on this
paper during the program on Semiclassical Methods at MSRI in 2003,
while the author  supported by the Clay Mathematics Institute.

\section{\label{HOL} Toeplitz quantization of symplectic maps}

In this section, we consider the quantization of symplectic maps
as Toeplitz Fourier integral operators. In some sense, the
material in this section is known, but it seems worthwhile to
recall the material and to complete some of the arguments.

 Suppose that $\chi: (M, \omega)
\to (M, \omega)$ is a symplectic diffeomorphism. There are several
equivalent ways to state the condition that $\chi$ is quantizable.
The most `geometric' one is the following:
\begin{defn} $\chi$ is quantizable if $\chi$ lifts to a contact transformation $\tilde{\chi}$
of $(X, \phi)$, i.e. a diffeomorphism of $X$ such that
$\tilde{\chi}^* \phi = \phi.$
\end{defn}
Equivalently, $\chi$ lifts to an automorphism of each power $L^N$
of the prequantum complex line bundle. It is said to be
linearizable in algebraic geometry.

Let us consider the obstruction to lifting a symplectic map. We
follow in part the discussion in   \cite{GS}, p. 220.  The key
notion is that $\chi$ preserve the holonomy map of the connection
$1$-form $\alpha$. Recall that the horizontal sub-bundle $H
\subset TX$  of the connection is defined by $H_x = \ker \alpha_x
= \{v \in T_x X: \alpha_x (v) = 0\}. $ The holonomy map
$$H : \Lambda \to U(1),\;\;\; H(\gamma) = e^{i \theta_{\gamma}} $$
from the free loop space   defined by horizontally lifting a loop
$\gamma: [0, 1] \to M$ to $\tilde{\gamma}: [0, 1] \to X$ and
expressing $\tilde{\gamma}(1) = e^{i \theta_{\gamma}}
\tilde{\gamma}(0). $ We say that $\chi$ is holonomy-preserving if
\begin{equation} H(\chi(\gamma)) = H(\gamma),\;\;\; \forall \gamma
\in \Lambda. \end{equation} If the loop is contained in the domain
of a local frame $s: U \to X$, then
\begin{equation} H(\gamma) = \exp (2 \pi i \int_{\gamma} s^* \alpha ).
\end{equation}
If $\gamma = \partial \sigma$, then $\int_{\gamma} s^* \alpha =
\int_{\sigma} \omega$. It follows (see \cite{GS}) that  symplectic
map preserves the holonomy around such homologically trivial
loops. Hence it is sufficient to consider the map
\begin{equation} H_{\chi} : H^1(M, \Z) \to U(1),\;\;\; H_{\chi} (\gamma) = H(\gamma)^{-1}  H(\chi(\gamma)) = e^{i
(\theta_{\gamma} - \theta{\chi(\gamma))}}  \end{equation}

\begin{prop} \label{SYMPLIFT} A symplectic map  $\chi$ of a symplectic manifold $(M, g)$ lifts to a contact transformation of
the associated prequantum $S^1$ bundle $(X, \alpha)$ if and only
if $H_{\chi} \equiv 1,$ the trivial representation. \end{prop}

\begin{proof} Suppose that $\chi$ lifts to $\tilde{\chi}: X \to
X$ as a contact transformation. Let $\gamma \in \Lambda$ and let
$\tilde{\gamma}$ be a horizontal lift of $\gamma$. Then
$\tilde{\chi}(\tilde{\gamma})$ is a horizontal lift of
$\chi(\gamma)$. Obviously, $\tilde{\gamma}(1) = e^{i
\theta_{\gamma}} \tilde{\gamma}(0) $ implies $\tilde{\chi} \circ
\tilde{\gamma}(1) = e^{i \theta_{\gamma}}\tilde{\chi} \circ
\tilde{\gamma}(0), $ so $H = H \circ \chi.$

Conversely, suppose that $H_{\chi} =1 $. We then define
$\tilde{\chi}$ by lifting $\chi$  along paths.  We fix a basepoint
$x_0 \in MX$ and define $\tilde{\chi}$ on the orbit $S^1 \cdot
x_0$ by fixing $\tilde{\chi}(x_0)$   to be a chosen basepoint on
$\pi^{-1}(\chi(\pi(x_0))$  and then extending by $S^1$ invariance.
We now consider horizontal paths $x(t): [0, 1] \to X$ from $x_0$.
At least one  horizontal path exists from $x_0$ to any given point
since the curvature is positive (Chow's theorem).  We define
$\tilde{\chi}(x(t))$ to be the horizontal lift of $\chi (\pi
x(t))$ to $\tilde{\chi}(x_0)$. To see that this is well-defined,
we must prove independence of the path. So let $x_1(t), x_2(t)$ be
two horizontal paths from $x_0$ to  $x_1(1)$. Thus, there is
trivial holonomy of the loop defined by $x_1$ followed by
$x_2^{-1}$ (i.e. the backwards path to $x_2$). Now project each
path, apply $\chi$, and horizontally lift. This defines a
horizontal lift of the  loop formed by the projected curves $\chi
\circ \pi \circ x_j(t)$ ($j = 1,2$). It  has trivial holonomy if
$\chi$ is holonomy preserving. It follows that the horizontal
lifts must agree at $t = 1$.

\end{proof}

It follows that $\chi$ always lifts to a contact transformation if
$M$ is simply connected. Hence,  if we lift $\chi$ first to the
universal cover $\tilde{M}$ of $M$, then this further  lifts to
$\tilde{X}$ as a contact transformation. We verify this in another
way, since we will use it in \S \ref{NOTSC}:

\begin{prop}\label{LIFTMAP} Let $\chi$ be a symplectic map of $\tilde{M}$.
Then  there exists a unique (up to one scalar) lift $\tilde{\chi}$
of $\chi$ to $\tilde{X}$ such that
$$ \pi \tilde{\chi} = \chi \pi $$
where $\pi: \tilde{X} \to \tilde{M}$ is the $S^1$-fibration.
\end{prop}

\begin{proof}

The fact that $\chi$  can be lifted to $\tilde{X}$ is obvious
since $\tilde{X} \simeq \tilde{M} \times S^1$. The key point is
that the map can be lifted as contact transformations. Any lift
which commutes with the $S^1$ action has the form
\begin{equation} \tilde{\chi}  \cdot (z, e^{i \theta}) = (\chi( z), e^{i
\theta + \phi_{\chi}(\theta, z) }). \end{equation}

The contact form on $\tilde{X}$ is the connection $1$-form
$\tilde{\alpha}$ of the hermitian line bundle over $\tilde{X}$. In
local symplectic coordinates $(x, \xi)$ on $\tilde{M}$ it has the
form
$$\tilde{\alpha} = \frac{1}{2} (\xi dx - x d\xi) - d \theta. $$
Since $\chi^*(x d \xi - \xi d x) - (x d \xi - \xi d x)$ is closed
on $\tilde{M}$, and since $\tilde{M}$ is simply connected, there
exists a function $f_{\chi} \in C^{\infty}(\tilde{M})$ such that
$$\tilde{\chi}^*(x d \xi - \xi d x) - (x d \xi - \xi d x) = d
f_{\chi}(x, \xi).$$ Using the product structure, we have that
$$\phi_{\chi}(x, \xi, \theta) = f_{\chi}(x, \xi)$$
defines a lift satisfying $\tilde{\chi}^*(\tilde{\alpha}) =
\tilde{\alpha},$ as desired.

Regarding uniqueness: the only flexibility in the lift is in the
choice of $f_{\chi}$, which is defined up to a constant. The
constant can be fixed by requiring that $f_{\chi}(0,0) = 0.$

\end{proof}

There is a weaker condition which has come up in some recent work
(cf. \cite{MR}): let us say that $\chi$ is quantizable at level
$N$ if $\chi$ lifts to an automorphism of the bundle $L^N$. Often
a map is quantizable of level $N$ along an arithmetic progression
$N = k N_0, k = 1,2,3, ...$ of powers although it is not
quantizable for all $N$. In geometric terms, this simply means
that $\chi$ fails to lift as a contact transformation of $X$ but
does lift as a contact transformation of $X/\Z_N$ where $\Z_N
\subset S^1$ is the group of $N$th roots of unity. In everything
that follows, the stated results have analogous for this modified
version of quantization.

\section{ \label{SC} Proof of Theorem \ref{DSMEETTOP} in the case $H^1(M, \C) = \emptyset$}

We first prove that if $H^1(M, \C) = \emptyset$, then every
automorphism is given by conjugation with a Toeplitz Fourier
integral operator. In this case, we may identify maps on $M$ which
$S^1$ invariant maps on $X$. We emphasize that we are not
considering the most general Toeplitz operators $\Pi A \Pi$ with
$A \in \Psi^*(X)$ but only the $S^1$-invariant operators whose
symbols lie in $C^{\infty}(M)$.  The proof is modelled on that of
Duistermaat-Singer \cite{DS1}, but has several new features due to
the holomorphic setting. In some respects the proof is simpler,
since there are no transmission automorphisms, and there are
natural identifications between symbols of different orders.
However in some respects it is more complicated, and also we must
be careful about using
 contravariant versus covariant symbols.

We begin with:

\begin{lem} \label{DS} Suppose that
$H^1(M, \C) = \emptyset$ and that  $\iota$ is an order-preserving
 automorphism of $ {\mathcal T}^{\infty} /{\mathcal T}^{-\infty}.$  Then  $\iota$ is equal to conjugation by a Toeplitz
 Fourier integral operator in the sense of Definition \ref{TFIO}.  \end{lem}

\begin{proof} Since $\iota$ is order preserving it induces automorphisms on
the quotients of the filtered algebra ${\mathcal T}^*. $

 We first consider ${\mathcal T}^0/ {\mathcal T}^{-1}$.  The map to contravariant
symbols defines an identification with $C^{\infty}(M).$ Thus
$\iota$ induces an automorphism of $C^{\infty}(M)$, viewed as an
algebra of contravariant symbols under multiplication. The maximal
ideal space of $C(M)$  equals $M$, hence $\iota$ induces a map
$\chi$ on $M$  such that $\iota(p) = p \circ \chi.$ Precisely as
in \cite{DS1} one verifies that $\chi$ is a smooth diffeomorphism
of $M$.

Now consider the quotients ${\mathcal T}^m/ {\mathcal T}^{m-1}.$ They are
are simply $N^m$ times   ${\mathcal T}^0/ {\mathcal T}^{-1}$, so for any $m$
 $\iota (p) = p \circ \chi$ for $p \in {\mathcal T}^m/ {\mathcal T}^{m-1}.$
This step is simpler than in the pseudodifferential case, and as a result certain
steps carried out in \cite{DS1} are unnecessary here.

Now let $n = 1$, so that ${\mathcal T}^1/ {\mathcal T}^{0}$ is a Lie algebra under
commutator bracket. The principal symbol is an isomorphism of the quotient algebra
to the Poisson algebra $(M, \{, \})$ defined by the symplectic form $\omega.$ Since
$\iota$ is an automorphism of the quotient algebra, we have
$$\{a \circ \chi, b \circ \chi\} = \{ a, b\} \circ \chi,$$
hence $\chi$ is a symplectic map of $(M, \omega).$ This step is also simpler than
in \cite{DS1}, and we see that no transmissions arise as possible automorphisms.

  By Proposition \ref{LIFTMAP},   $\chi$ lifts to a contact
transformation $\tilde{\chi}$ of $X$.

\subsubsection{\label{SPA} Symbol preserving automorphisms}

Now let $A_N^{-1} = \Pi_N a T_{\chi}^{-1} \Pi_N$ denote any
Toeplitz quantization of $\tilde{\chi}^{-1}.$ It follows that
$\alpha(P) = A_N \iota(P) A_N^{-1}$ is an automorphism of
${\mathcal T}^{\infty}/ {\mathcal T}^{- \infty}$ which preserves
principal contravariant symbols in the sense of \S \ref{COCO}. We
now prove that any such automorphism is given by conjugation with
a Toeplitz multiplier of some order $s$.

Thus, let $j$ be a principal contravariant symbol preserving automorphism.
Let $P \in {\mathcal T}^m$.  Since $j(P) - P \in {\mathcal T}^{m - 1}$, we get
an induced map
\begin{equation} \beta_m : C^{\infty}(M) \to C^{\infty}(M),\;\;\;\ \beta_m(a) = j(\Pi_N a \Pi_N)  - \Pi_N a \Pi_N. \end{equation}
Then $\beta = \beta_m$  is a derivation in two ways:
\begin{equation} \begin{array}{ll}(i) &  \beta(p \cdot q) = \beta(p) \cdot q + p \cdot \beta(q)\\ & \\(ii) &  \beta(\{p,q\}) = \{\beta(p), q\} + \{p, \beta(q)\}.\end{array}
\end{equation}

Now any derivation in the sense of $(i)$ is given by differentiation along a vector
field $V.$ Since $V$ commutes with Poisson bracket, it must be a symplectic vector
field. Since $\omega(V, \cdot) =
\beta$ is a closed $1$-form, there exists a local Hamiltonian $H$ for $V$. Under
our assumption that $H^1(M) = \{0\}$, the Hamiltonian is global, so $V = \Xi_H$ for some global $H$. Thus we have:
\begin{equation} \beta_m(a) = i \{a, H_{\log b}\} \end{equation}
for some $b_m \in C^{\infty}(M).$ Thus, $j$ may be represented by
$j(P) = B^{-1} P B$ for $P \in {\mathcal T}^m$, with $B = \Pi_N
e^{i b} \Pi_N \; (b = b_m).$ Because of the natural identification
of $S_{scl}^m \equiv N^m S_{scl}^0$, we find that $b_m = b$ is the
same for all $m$.

By composing automorphisms, we now have an automorphism $j_2$ such that
$$j_2(\Pi_N a_N \Pi_N) - \Pi_N a_N \Pi_N \in {\mathcal T}^{m - 2},\;\;\; a_N \in S_{scl}^m.$$
We find as above that $j_2(\Pi_N a_N \Pi_N) - \Pi_N a_N \Pi_N = \beta_2(a_N)$ where
$\beta_2$ is a derivation, hence $\beta_2 = N^{-1} \{\log b_{-1}, \cdot\}$ for some $b_{-1}$.
Therefore,
\begin{equation} \Pi_N e^{- i N^{-1} b_{-1} } \Pi_N j_2(P) \Pi_N e^{ i N^{-1} b_{-1} }
\Pi_N
- P \in {\mathcal T}^{m - 2},\;\;\; \forall P \in {\mathcal T}^m. \end{equation}

Proceeding in this way, we get an element $b_N \in S_{scl}^0$ such that
\begin{equation} \Pi_N e^{- i b_N } \Pi_N j(P) \Pi_N e^{ i b_N }\Pi_N
- P \in {\mathcal T}^{-\infty},\;\;\; \forall P \in {\mathcal T}^m. \end{equation}
 This completes
the proof of the Lemma.

\end{proof}

\subsubsection{Conclusion of proof when $H^1(M, \C) = \emptyset$}

Lemma 3 of \cite{DS1} is an abstract result which says that automorphisms of Frechet
spaces satisfying a certain density condition are always given by conjugation. The
density condition is easy to prove, so we omit the  proof. The result is:
\medskip

{\it Let $\iota$ denote an automorphism of ${\mathcal T}^*$ acting on $H^{\infty}(M).$
Then $\iota(P) = A^{-1} P A,$ where $A: H^{\infty}(M) \to H^{\infty}(M)$ is an invertible,
continuous linear map, determined uniquely up to multiplicative constant. }
\medskip

Thus,  $i(P) = A^{-1} P A$
for all $P \in {\mathcal T}^*$ and also, by Lemma 1, there exists an elliptic Toeplitz Fourier
integral operator $B$ such that $i(P) \equiv B^{-1} \circ P \circ B$ for all $P \in
{\mathcal T}^{\infty}/{\mathcal T}^{- \infty}.$
Let $E = A \circ B^{-1}$. Then $[E, P]  \in {\mathcal T}^{-\infty}$ for all $P
\in {\mathcal T}^{\infty}/{\mathcal T}^{- \infty}.$ In particular, $[E, P] \in
{\mathcal T}^{-\infty}$
for all $P = \{\Pi_N a \Pi_N \}, a \in C^{\infty}(M).$

In place of \cite{DS1}, Lemma 4, we use
\begin{lem} \label{E} Let $E$ be an operator
on $ {\mathcal H}^2$ such that $[E, \Pi a \Pi] \in {\mathcal T}^{-
\infty}$ for all $a \in C^{\infty}(M)$. Then there exists a
constant $c$ such that  $E = c \Pi + R$, where $R$ is a smoothing
operator. \end{lem}

\begin{proof} It is sufficient to prove the statement for all $a$
supported in a given $S^1$-invariant open set $U \subset X$. We
can then use a partition of unity to prove the result for all $a$.
We use the notation $A \sim_U B$ to mean that $A, B$ are defined
on $U$ and their difference is a smoothing operator on $U$.

 In a
sufficiently small open set $U \subset X$, there exists a Fourier
integral operator $F: L^2(X) \to L^2(\R^n)$ associated to a
contact transformation $\phi$  such that $ \Pi \sim_U F \Pi_0 F^*$
modulo smoothing operators, where $\Pi_0$ is the model \szego
kernel discussed in \cite{BS, BG}, namely the orthogonal
projection  onto the kernel of the annihilation operators $D_j =
\frac{1}{i} (\partial/\partial y_j + y_j |D_t|)$ on $\R^n = \R^p
\times \R^q$. Furthermore, it is proved in \cite{BS, BG} that
there exists a complex Fourier integral operator $R_0: L^2(\R^q)
\to L^2(\R^n)$ such that $R_0^* R_0 = I, R_0 \; R_0^* = \Pi_0.$
Moreover for any pseudodifferential operator $A$ on $X$, there
exists a pseudodifferential operator $Q$ on $\R^q$ so that $\Pi_0
A \Pi_0 \sim R_0 Q R_0^*.$ Transporting $R_0$ to $X$ by $F$, we
obtain a complex Fourier integral operator $R : L^2(\R^q) \to
L^2(X)$  so that $R R^* \sim_U \Pi, R^* R \sim_{\phi(U)} I$  and
so that $\Pi a \Pi \sim_U R Q R^*$. Then $[E, \Pi a \Pi] \sim_U
[E, R^* Q R]$ and we may rewrite the condition on $E$ as:
\begin{equation} [E, R^* Q R] \in {\mathcal T}^{-\infty}(X),\;\;
\forall Q \in \Psi^0(\R^q). \end{equation} This is equivalent to
\begin{equation} [R^* E R, Q ] \in \Psi^{-\infty}(X),\;\;
\forall Q \in \Psi^0(\R^q). \end{equation} We then apply Beals'
characterization of pseuodifferential operators:
\medskip
{\it $P \in \Psi^k (\R^m)$ if and only if for all $\{j_i,
k_{\ell}\},$
$$ad(x_{j_1}) \cdots ad(x_{j_r}) ad (D_{x_{k_1}}) \cdots ad (D_{x_{k_s}}) P\; H^{s + r}(\R^m) \to H^{s }(\R^m)$$
is bounded.} Here,  $ad(L) P$ denotes $  [L, P].$ It follows first
that  $R^* E R \in \Psi^0$, and easy symbol  calculus shows that
the complete symbol of $R^* E R$ is  constant. Hence, $R^* E R = I
+ S$ where $S$ is a smoothing operator. Applying $R$ on the left
and $R^*$ on the right concludes the proof.

\end{proof}

\begin{rem} It is in the step in \S \ref{SPA} that the distinction
between symplectic maps of $M$ and contact transformations of $X$
enters. Ultimately it is this step which leads to Corollary
\ref{PHYSICS}. We also note that the proof above is rather
different from that in \cite{DS1}.
 \end{rem}

\section{\label{NOTSC} $H^1(M, \C) \not= \emptyset$}

The problems with quantizing symplectic maps on $M$  are all due
to the fundamental group $\pi_1(M)$ or more precisely $H_1(M,
\C)$.  We solve them by passing to the universal cover
$\tilde{M}$. In this section, we relate Toeplitz operators on $M$
and $\tilde{M}$.

\subsection{Toeplitz operators  on the universal cover}

Since we are comparing algebras and automorphisms on covers to
those on a quotient, we begin with the abstract picture as
discussed in \cite{GHS}.  We then specialize it to algebras of
Toeplitz operators.

\subsubsection{Abstract theory}  Suppose that $p: \tilde{X} \to X$
is a covering map of a compact manifold $X$, and denote its deck
transformation group by $\Gamma$.  We regard $\Gamma$ as acting on
the left of $\tilde{X}$.  We would like to compare operators on
$\tilde{X}$ and operators on $X$. To gain perspective, we start
with the large von Neumann algebra ${\mathcal B}_{\Gamma}$ of all
bounded operators on ${\mathcal L}^2(\tilde{X})$ which commute
with $\Gamma$. Later we specialize to the Toeplitz algebra which
is our algebra of observables.

The Schwartz kernel of such an operator satisfies $B(\gamma x,
\gamma y) = B(x, y)$. If we denote by ${\mathcal D}$ a fundamental
domain for $\Gamma$, then there exists an identification
\begin{equation} \label{ISO} {\mathcal L}^2(\tilde{X}) \simeq {\mathcal L}^2(\Gamma)
\otimes {\mathcal L}^2({\mathcal D}) = {\mathcal L}^2(\Gamma)
\otimes {\mathcal L}^2(X).
\end{equation} Elements of ${\mathcal L}^2(\Gamma) \otimes {\mathcal L}^2(X)$ can be
viewed as functions $f(\gamma, x) $ on $\Gamma \times {\mathcal
D}$.  The unitary isomorphism is defined by
$$\phi \in {\mathcal L}^2(\tilde{X}) \to f_{\phi}(\gamma, x) = \phi(\gamma \cdot x). $$
Note that both left translation $L_{\gamma}$ and right translation
$R_{\gamma}$ by $\gamma$ act on this space, namely $$L_{\gamma}
f(\alpha, x) = f(\gamma \alpha, x),\;\;\; R_{\gamma} f(\alpha, x)
= f(\alpha \gamma, x). $$ We may regard ${\mathcal B}_{\Gamma}$ as
bounded operators commuting with all $L_{\gamma}$. The isomorphism
(\ref{ISO}) induces an algebra isomorphism
\begin{equation} \label{ALGISO} {\mathcal B}_{\Gamma} \simeq
{\mathcal R} \otimes {\mathcal B}(X),
\end{equation} where ${\mathcal B}(X)$ is the algebra of bounded operators on
$X$ and where ${\mathcal R}$ is the algebra generated by right
translations $R_{\gamma}$ on $L^2(\Gamma)$.

So far we have been considering operators on ${\mathcal
L}^2(\tilde{X})$. The corresponding  algebra ${\mathcal
B}_{\Gamma}$ is much larger than ${\mathcal B}(X)$.  To make the
connection to ${\mathcal L}^2(X)$ tighter, we
 need to consider the space ${\mathcal L}^2_{\Gamma}(\tilde{X})$ of
$\Gamma$-periodic functions on $\tilde{X}$. The natural Hilbert
space structure is to define $||f||_{\Gamma}^2  = \int_{{\mathcal
D}} |f(x)|^2 dV$ where $dV$ is a $\Gamma$-invariant volume form.
We have the obvious isomorphism ${\mathcal
L}^2_{\Gamma}(\tilde{X}) \simeq {\mathcal L}^2(X)$. We may regard
elements of ${\mathcal L}^2_{\Gamma}(\tilde{X}) $ as functions
$f(\gamma, x)$ as above which are constant in $\gamma$.

Elements   $B \in {\mathcal B}_{\Gamma}$ with properly supported
kernels, or kernels which decay fast enough off the diagonal,  act
on  ${\mathcal L}^2_{\Gamma}(\tilde{X}) $. Indeed,  ${\mathcal R}$
acts trivially on ${\mathcal L}^2_{\Gamma}(\tilde{X})$, so $
{\mathcal B}_{\Gamma}$ acts by the quotient algebra $ {\mathcal
B}_{\Gamma}/{\mathcal R}_{\Gamma}$. We will be working with
subalgebras of Toeplitz operators where the action is clearly
well-defined.

\begin{rem} Let us define ${\mathcal B}_{\Gamma}^{\Gamma}$ as the subalgebra
of ${\mathcal B}_{\Gamma}$ of elements which  commute with both
$R_{\gamma}$ and $L_{\gamma}$ for all $\gamma$. Then we have:
${\mathcal B}_{\Gamma}^{\Gamma} (\tilde{X}) \simeq {\mathcal
B}(X).$ We may write the (Schwartz) kernel of an element of
${\mathcal B}$ as $B(\gamma, x, \gamma', x'). $ It belongs to
${\mathcal B}_{\Gamma}$ if $B(\alpha \gamma, x, \alpha \gamma',
x') = B(\gamma, x, \gamma', x')$ and it belongs to ${\mathcal
B}_{\Gamma}^{\Gamma}$ if addditionally $B( \gamma \alpha, x,
\alpha, \gamma' \alpha, x') = B(\gamma, x, \gamma', x')$
\end{rem}

We have been talking about  algebras of bounded operators, but our
main interest is in $C^*$ algebras of Toeplitz operators.
Everything we have said restricts to these subalgebras  once we
have defined the appropriate notions.

\subsection{Toeplitz operators}

The positive hermitian holomorphic line bundle $(L, h) \to M$
pulls back under $\pi$ to one $(\tilde{L}, \tilde{h}) \to
\tilde{M}$. This induces an inner product on the space
$H^0(\tilde{M}, \tilde{L}^N)$ of entire holomorphic sections of
$\tilde{L}^N$. We denote by $ {\mathcal H}^2(\tilde{M},
\tilde{L}^N) $ the space of $L^2$ holomorphic sections relative to
this inner product.

As in the quotient, there exists an associated $S^1$ bundle
$\tilde{X}$ with a contact (connection) form $\tilde{\alpha}$ such
that
$$\begin{array}{lll} \tilde{X} & \rightarrow & \tilde{M} \\ & & \\
\downarrow & & \downarrow \\ & & \\
X & \rightarrow & M \end{array} $$ commutes. The vertical arrows
are covering maps and the horizontal ones are $S^1$ bundles. We
denote the deck transformation group of $\tilde{X} \to X$ by
$\tilde{\Gamma}$. It is isomorphic to $\Gamma$, so when no
confusion is possible we drop the $\tilde.$ Since all objects are
lifted from quotients, it is clear that $\tilde{\Gamma}$ acts by
contact transformations of $\tilde{\alpha}$. Let us the denote
operator of translation by $\gamma$ on $\tilde{M}$ by
$L_{\gamma}.$

We denote by ${\mathcal H}^2(\tilde{X})$ the Hardy space of $L^2$
CR functions on $\tilde{X}$. They are boundary values of
holomorphic functions in the strongly pseudoconvex complex
manifold
\begin{equation} \tilde{D}^* =\{(z, v) \in \tilde{L}^*: h_z(v)< 1
\} \end{equation} which are ${\mathcal L}^2(\tilde{X})$. The group
$\tilde{\Gamma} $ acts on  $\tilde{D}^* \subset \tilde{L}^*$ with
quotient the compact  disc bundle $D^* \subset L^* \to M$. In this
setting it is known (cf. \cite{GHS}, Theorem 0.2) that
$$\dim_{\Gamma}  {\mathcal H}^2(\tilde{X})  = \infty. $$

Due to the $S^1$ symmetry, holomorphic functions on $\tilde{D}^*$
are easily related to CR holomorphic functions on $\tilde{X}$. We
denote by
\begin{equation}\tilde{\Pi}: {\mathcal L}^2(\tilde{X}) \to {\mathcal H}^2(\tilde{X}) \end{equation}
the \szego (orthogonal) projection.  Under the $S^1$ action,  we
have
\begin{equation} {\mathcal H}^2(\tilde{X}) = \bigoplus_{N=1}^{\infty} {\mathcal
H}_N^2(\tilde{X}),\;\;\;\tilde{\Pi} =  \bigoplus_{N=1}^{\infty}
\tilde{\Pi}_N.
\end{equation}
As on $X$, we have ${\mathcal H}_N^2(\tilde{X}) \simeq {\mathcal
H}^2(\tilde{M}, \tilde{L}^N).$ We refer to \cite{GHS} (see example
2) on p. 559)  for the proof that $\tilde{L} \to \tilde{M}$ has
many holomorphic sections.

So far we have discussed $L^2$ functions on $\tilde{X}$. More
important are periodic functions. We endow them with a Hilbert
space structure by setting:\begin{equation} \left\{
\begin{array}{l} {\mathcal L}^2_{\Gamma}(\tilde{X}) = \{f \in L^2_{loc} (\tilde{X}): L_{\gamma} f =
f\}, \\ \\
{\mathcal H}^2_{\Gamma}(\tilde{X}) = \{f \in {\mathcal
L}^2_{\Gamma}(\tilde{X}), \overline{\partial}_b f = 0 \},
\end{array} \right.
\end{equation} with the inner product $\langle , \rangle_{\Gamma}$
obtained by integrating over a fundamental domain ${\mathcal D}$
for $\Gamma$. Both are  direct sums of weight spaces ${\mathcal
L}^2_{\Gamma, N}, $ resp. $ {\mathcal H}^2_{\Gamma, N}(\tilde{X})$
for the $S^1$ action on $\tilde{X}$. There  exists a Hilbert space
isomorphism $L:
 {\mathcal H}^2_{\Gamma}(\tilde{X}) = {\mathcal H}^2(X)$, namely
by lifting $L \phi = p^* \phi$ under the covering map. The adjoint
of $L$  is given by:
$$L^* f = p_* (f 1_{\mathcal D}). $$
We thus have
\begin{equation} L L^* = Id:  {\mathcal H}^2_{\Gamma}(\tilde{X})
\to  {\mathcal H}^2_{\Gamma}(\tilde{X}). \end{equation} We observe
that the spaces ${\mathcal H}^2_{\Gamma, N}(\tilde{X})$ and
${\mathcal H}^2(\tilde{X})$ are completely unrelated and have
different dimensions.

 We now consider Toeplitz algebras. Since $\tilde{X}$ is
 non-compact in general, we must take some care that Toeplitz
 operators are well-defined. Otherwise, the definitions are the
 same as for $X$: ${\mathcal T}^s(\tilde{X})$ is the space of
 operators of the form $\tilde{\Pi} A \tilde{\Pi}$ where $A \in
 \Psi_{S^1}^s(\tilde{X})$ is the space of properly supported
 pseudodifferential operators commuting with $S^1$. We also define
 ${\mathcal T}^{-\infty}$ as the space of such $\tilde{\Pi} A \tilde{\Pi}$ with
 $A$ having a smooth properly supported kernel.

 We then distinguish the automorphic Toeplitz operators:
$${\mathcal T}_{\Gamma}^*(\tilde{X}) = \{\tilde{\Pi} A \tilde{\Pi} \in {\mathcal T}^*(\tilde{X}):
L_{\gamma}^*\tilde{\Pi} A \tilde{\Pi} L_{\gamma} = \tilde{\Pi} A
\tilde{\Pi} \}.$$ We note that $[L_{\gamma}, \tilde{\Pi}] = 0$ or
equivalently $\tilde{\Pi} (\gamma x, \gamma y) = \tilde{\Pi}(x,
y).$ So the operative condition is that $A \in \Psi_{S^1,
\Gamma}$, the space of pseudodifferential operators commuting with
$\Gamma$. The associated symbols $a_N(z, \bar{z})$ are exactly the
$\Gamma$-invariant symbols on $\tilde{M}$. We note:

\begin{prop}\label{ACT}  There exists a  representation $\rho_{\Gamma}$ of ${\mathcal T}_{\Gamma}^*(\tilde{X})$
on ${\mathcal H}^2_{\Gamma}$ preserving each ${\mathcal H}
^2_{\Gamma, N}$.
\end{prop}

\begin{proof} Since  $L_{\gamma} \Pi A \Pi f = \Pi A \Pi f$
whenever $L_{\gamma} f = f$, the only issue is whether $A f$ is
well-defined for $f \in {\mathcal H}_{\Gamma}$. However, $A$ is a
polyhomgeneous sum of $a_j {\mathcal N}^{-j}$ where $a_j $ is a
periodic function, so the action is certainly defined.

\end{proof}

Let us denote by ${\mathcal K}_{\Gamma} = \ker \rho_{\Gamma}$. We
further denote by $\rho_{\Gamma, N}$ the associated representation
on ${\mathcal H}_{\Gamma, N}^2$, and put ${\mathcal K}_{\Gamma} =
\ker \rho_{\Gamma, N}$.

\begin{prop} ${\mathcal K}_{\Gamma} \subset {\mathcal
T}^{-\infty}(\tilde{X})$.
\end{prop}

\begin{proof}

Assume  that $\Pi A \Pi \in {\mathcal T}_{\Gamma}^*(\tilde{X})$
annihilates ${\mathcal H}^2_{\Gamma}(\tilde{X}) \simeq {\mathcal
H} ^2(X)$. This means that $\langle A f, g \rangle_{{\mathcal D}}
= 0$ for all $f, g \in {\mathcal H}^2_{\Gamma}(\tilde{X}).$ In
particular it implies that the `Berezin symbol' $\tilde{\Pi}_N a_N
\tilde{\Pi}_N (z,z) = 0$. However, asymptotically, $N^{-m}
\tilde{\Pi}_N a_N \tilde{\Pi}_N (z,z) \sim a_0(z)$ for a zeroth
order Toeplitz operator. One sees by induction on the terms in
(\ref{SPECIAL}) that $A \sim 0$.

\end{proof}

It could  happen that  ${\mathcal K}_{\Gamma}  \not= 0$, unlike
the analogous representation on ${\mathcal H}^2(X)$ which defines
Toeplitz operators.   For each  $N$ there could exist $a_N \in
C^{\infty}(M)$ with $||a_N||_{L^2} = 1$ which is orthogonal to the
finite dimensional space ${\mathcal H}_{\Gamma, N }^2(\tilde{X})
\simeq {\mathcal H}^2_N (X)$ but which is not orthogonal to
${\mathcal H}^2_N(\tilde{X})$. Then  $T = \tilde{\Pi}_N a_N
\tilde{\Pi}_N \in {\mathcal K}_{\Gamma, N}.$

 We now relate ${\mathcal T}_{\Gamma}^*(\tilde{X}) $ to
${\mathcal T}^*(X)$. In preparation, we relate the \szego kernels
on $X, \tilde{X}$. First, we consider a fixed $N$. The following
is proved in \cite{SZ}:

\begin{prop} \label{UPSZEGON} The degree $N$ \szego kernels of $X, \tilde{X}$ are related by:
$$\Pi_N(x, y) = \sum_{\gamma \in \Gamma} \tilde{\Pi}_N(\gamma \cdot
x, y). $$ The same formula defines the \szego projector
$L^2_{\Gamma, N} \to {\mathcal H}^2_{\Gamma, N}$.
\end{prop}

The key point is to use the estimate
\begin{equation} \label{OFFDIAG}  |\tilde{\Pi}_N(
x, y)| \leq C e^{-\sqrt{N} \tilde{d}(x, y)} \end{equation} where
$\tilde{d}(x,y)$ is the Riemannian distance with respect to the
\kahler metric $\tilde{\omega}$ to show that the  sum converges
for sufficiently large $N$. The estimates show that
$\tilde{\Pi}_N$ acts on ${\mathcal L}^{\infty}(\tilde{X})$. Since
${\mathcal L}^2_{\Gamma}(\tilde{X}) \cap C(\tilde{X}) \subset
{\mathcal L}^{\infty}(\tilde{X}), $ $\tilde{\Pi}_N$ acts on
$H^2_{\Gamma}(\tilde{X})$. The formula is an immediate consequence
of writing $$\tilde{\Pi}_N s_N(x) = \int_{\tilde{X}}
\tilde{\Pi}_N(x, y) s_N(y) dV(y) = \int_{{\mathcal D}}
\sum_{\gamma \in \Gamma} \tilde{\Pi}_N(x, \gamma y) s_N(y) dV(y).
$$ for a periodic $s_N$ in terms of a fundamental domain.

Now we consider the full \szego kernel:

\begin{cor} \label{UPSZEGO} We have  $L \Pi_N = \tilde{\Pi}_N
L$ for each $N$ as operators from ${\mathcal L}^2_N(X)$ to
${\mathcal H}^2_{\Gamma, N}$.  Hence, $L \Pi= \tilde{\Pi} L:
{\mathcal L}^2(X) \to {\mathcal H}^2_{\Gamma}(\tilde{X})$.

\end{cor}

 The identity (\ref{ALGISO}) for the larger algebra of
bounded operators suggests that ${\mathcal
T}_{\Gamma}^*(\tilde{X})$ should be a larger algebra than
${\mathcal T}^*(X).$  However, this is not the case.

\begin{prop} $L$ induces an algebra isomorphism ${\mathcal
T}^*(X) \simeq {\mathcal T}_{\Gamma}(\tilde{X})/{\mathcal
K}_{\Gamma}(\tilde{X}).$
\end{prop}

\begin{proof} It follows from Proposition \ref{ACT} and Corollary \ref{UPSZEGO}  that
\begin{equation} \label{TOEPISO} L \Pi a \Pi = \tilde{\Pi} p^*a
\tilde{\Pi} L : L^2(X) \to {\mathcal H}^2_{\Gamma}(\tilde{X}) \iff
\Pi a \Pi = L^* \tilde{\Pi} p^*a \tilde{\Pi} L : {\mathcal L}^2(X)
\to {\mathcal H}^2_{\Gamma} (X) .\end{equation}

Equality on the designated spaces is equivalent to equality in the
algebras. Further,  the equality $L L^* = I: {\mathcal
H}^2_{\Gamma} \to {\mathcal H}^2_{\Gamma}$ implies that the linear
isomorphism is an algebra isomorphism.

\end{proof}

It may seem surprising that ${\mathcal T}_{\Gamma}(\tilde{X})$ is
so `small'.  We recall that one  obtains ${\mathcal T}^*(X)$  by
representing  ${\mathcal T}_{\Gamma}(\tilde{X})$  on ${\mathcal
H}^2_{\Gamma}(\tilde{X})$. When dealing with all bounded
operators, the kernel is very large (${\mathcal R}_{\Gamma}$) and
it is also large if we fix $N$ and consider the associated
Toeplitz algebra. But the kernal is trivial if we consider the
full Toeplitz algebra.

We also have:

\begin{cor} $L$ induces algebra isomorphisms ${\mathcal T}_N(X)
\simeq {\mathcal T}_{\Gamma, N}(\tilde{X})/{\mathcal K}_{\Gamma,
N}.$
\end{cor}

Now let us relate such automorphisms to automorphisms on
$\tilde{X}$.   The concrete identification with ${\mathcal T}(X)$
is by (\ref{TOEPISO}). We have:

\begin{prop}\label{AUTO}  There is a natural identification of:
\begin{itemize}

 \item automorphisms $\alpha$ on ${\mathcal T}^*(X)$ with
 automorphisms $\tilde{\alpha}$ on ${\mathcal T}^{* }_{\Gamma}(\tilde{X})/{\mathcal
 K}_{\Gamma}$:

\item automorphisms $\alpha$ on ${\mathcal T}^*(X)/ {\mathcal
T}^{-\infty}(X)$ with
 automorphisms $\tilde{\alpha}$ on ${\mathcal T}^{* }_{\Gamma}(\tilde{X})/{\mathcal
 T}^{-\infty}_{\Gamma}$.

 \end{itemize}
 \end{prop}

 \begin{proof} The first statement is clear since the  algebras are the
 isomorphic. An  automorphism of
 ${\mathcal T}^*_{\Gamma}$ descends to ${\mathcal T}^*$ if and
 only if it preserves ${\mathcal K}_{\Gamma}$. Concretely, we
 wish to set:
\begin{equation} \label{CONCISO} \alpha(\Pi a \Pi) = L^* \tilde{\alpha} (\tilde{\Pi}  p^*a
 \tilde{\Pi}) L. \end{equation}
 The inner operator must be determined by the left side for this
 to be well-defined. We  have
 \begin{equation} \tilde{\Pi}  p^*a
 \tilde{\Pi} \equiv L \Pi a \Pi L^* \;\; \mbox{mod}\;\; {\mathcal K}_{\Gamma}. \end{equation}
 The same equivalence is true modulo the larger subalgebra
 ${\mathcal T}_{\Gamma}^{-\infty}.$
Since ${\mathcal
 T}^{-\infty}_{\Gamma}/K_{\Gamma} ={\mathcal T}^{-\infty} (X)$, the second
 statement is correct.

 \end{proof}

\subsection{Completion of the proof of Theorem \ref{DSMEETTOP} }

We now prove statements (iii)-(iv) of the Theorem.  The following
gives an `upper bound' on existence of semi-classical
automorphisms.

\begin{lem} \label{NONSIMAUTO} Suppose that $\alpha$ is an order-preserving automorphism
 of ${\mathcal T}^*(X)/{\mathcal T}^{-\infty}$.
Let $\tilde{\alpha}$ denote the corresponding automorphisms of
${\mathcal T}^{*}_{\Gamma}(\tilde{X})/{\mathcal K}_{\Gamma}$.
Then there exists a canonical transformation $\chi$ of $(M,
\omega)$ and a unitary Toeplitz quantum map $\tilde{U}_{\chi}$ on
$H^2(\tilde{X})$ such that $$ \tilde{\alpha} (P) =
\tilde{U}_{\chi}^* P \tilde{U}_{\chi}, $$ and such that
$$L_{\gamma}^{-1} \tilde{U}_{\chi} L_{\gamma} = M_{\gamma}
\tilde{U}_{\chi},$$ where $M_{\gamma} \in {\mathcal T}
^*_{\Gamma}(\tilde{M})'.$
\end{lem}

\begin{proof}
By Lemma(\ref{DS}), we know that $\tilde{\alpha}$ is given by
conjugation by a Toeplitz quantum map $\tilde{U}_{\chi}$. We may
define $M_{\gamma}$ by the formula above since $\tilde{U}_{\chi}$
is invertible. We then determine its properties.

We have: $$L_{\gamma}^{-1} \tilde{U}_{\chi}^{-1} L_{\gamma} P
L_{\gamma}^{-1}  \tilde{U}_{\chi} L_{\gamma} =
\tilde{U}_{\chi}^{-1} P  \tilde{U}_{\chi}, \;\; \forall P \in
{\mathcal T}^*_{\Gamma}(\tilde{M}).$$ Hence,
$$ M_{\gamma}^{-1} \tilde{U}_{\chi}^{-1}  P   \tilde{U}_{\chi} M_{\gamma} = \tilde{U}_{\chi}^{-1} P  \tilde{U}_{\chi} \iff M_{\gamma}^{-1} P M_{\gamma} = P  \;\;
\forall P \in {\mathcal T}^*_{\Gamma}(\tilde{M}) .$$ It follows
that $M_{\gamma}$ is central. This proves (iii) of Theorem
\ref{DSMEETTOP}.
\end{proof}

The `lower bound' is given by:

\begin{lem} Suppose that $\chi$ is a symplectic map of $(M,
\omega)$. Then it lifts to a contact transformation $\tilde{\chi}$
of $(\tilde{X}, \tilde{\alpha})$. The associated quantum map
$U_{\tilde{\chi}}$ on $\tilde{X}$ defines (by conjugation) an
order preserving automorphism $\tilde{\alpha}$ of ${\mathcal
T}_{\Gamma}(\tilde{X})$ which descends to  an automorphism of
${\mathcal T}^*(X)/{\mathcal T}^{-\infty}$. If $\tilde{\alpha}$
preserves ${\mathcal K}_{\Gamma}$, then it defines an automorphism
of all of ${\mathcal T}^*(X)$.
\end{lem}

\begin{proof} $\chi$ automatically lifts to $\tilde{M}$ as a
symplectic map commuting with the action of $\Gamma$. By
Proposition \ref{LIFTMAP} it lifts to $\tilde{X}$ as a contact
transformation. We then define a unitary quantum map by
$$\tilde{U}_{N \tilde{\chi}} = \tilde{\Pi}_N \sigma
T_{\tilde{\chi}} \tilde{\Pi}_N, $$ where $\sigma$ is a function on
$\tilde{M}$ which makes the operator unitary on ${\mathcal H}
^2(\tilde{X})$. See \cite{Zel1} for background.

A crucial issue now is the commutation relations between
$\tilde{U}_{N \tilde{\chi}} $ and $\tilde{\Gamma}$. When $\chi$
lifts to $X$, i.e. is quantizable,  then $\tilde{\chi}$ commutes
with $\tilde{\Gamma}$. In this case, $\chi$ is quantizable as a
quantum map and there was no need to lift it to $\tilde{X}$ to
quantize it as an automorphism.

Assume however that $\chi$ does not lift to $X$ and consider the
commutation relations of the translation by $\tilde{\chi}$ with
left translations by elements of $\tilde{\Gamma}$. The commutator
$\tilde{\chi} \tilde{\gamma} \tilde{\chi}^{-1}
\tilde{\gamma}^{-1}$ covers the identity map of $\tilde{M}$ since
the lift of $\chi$ to $\tilde{M}$ commutes with $\Gamma$.
Furthermore, it commutes with the $S^1$ action on $\tilde{X}$. It
follows that
\begin{equation} \tilde{\chi} \tilde{\gamma} \tilde{\chi}^{-1}
\tilde{\gamma}^{-1}  = \; T_{e^{i \theta_{\gamma, \chi}}}
\end{equation}
where the right side is translation by the element $e^{i
\theta_{\gamma, \chi}}$. The angle $\theta_{\gamma, \chi}$ is
apriori a function on $\tilde{M}$. However, the left side is a
contact transformation covering the identity  and therefore $d
\theta = d \theta + d \theta_{\gamma, \chi}$, i.e.
$\theta_{\gamma, \chi}$ is a constant.

After quantizing, the same commutator identity holds for the
operators. Therefore the operators $M_{\gamma}$ are central. It
follows that the   automorphism
\begin{equation} \tilde{\alpha}(P) = \tilde{U}_{\chi}^* P
\tilde{U}_{\chi}\end{equation}  satisfies $$\tilde{\alpha}_N
(\tilde{\Pi}_N p^*a_N
 \tilde{\Pi}_N )\in {\mathcal T}_{\Gamma}^*
 (\tilde{X}).$$
 It therefore descends to ${\mathcal T}^*(X)/{\mathcal T}^{-\infty}(X)$  by
(\ref{CONCISO}), i.e. as
$$\tilde{\alpha}_N
(\tilde{\Pi}_N p^*a_N
 \tilde{\Pi}_N ) = L^* \tilde{U}_{\chi}^* \tilde{\Pi}_N p^*a_N
 \tilde{\Pi}_N ) \tilde{U}_{\chi} L.$$
 Unitarity of $L$ then implies that
$\alpha$ is also an
 automorphism. If the automorphism preserves ${\mathcal
 K}_{\Gamma}$ then it also descends to ${\mathcal T}^*(X)$.

\end{proof}

An obvious question is whether the condition that the automorphism
preserve ${\mathcal K}_{\Gamma}$ is equivalent to the quantization
condition that $\chi$ lift to $X$. Clearly, quantizability in the
sense of Definition \ref{TQM}  implies preservation of ${\mathcal
K}_{\Gamma}$, since the quantum map $U_{\chi, N}$ is well defined
on the spaces ${\mathcal H}_N$. The converse is not obvious, since
we only know apriori that the automorphism induces automorphisms
$\alpha_N$ of the finite rank observables $Op_N(a_N)$ for fixed
$N$. Abstractly, such automorphisms must be given by conjugations
by unitary operators on ${\mathcal H}_N$, but it is not clear that
these unitary operators  are Toeplitz quantum maps in the sense of
\ref{TQM}.

We end the section with

\subsubsection{Proof of Corollary \ref{PHYSICS}}

\begin{proof} This follows immediately from (ii) if $H^1(M, \C) =
\emptyset.$ If $H^1(M, \C) \not= \emptyset,$ then by  Theorem
\ref{DSMEETTOP}  (iii) there exists a symplectic map $\chi$ of $M$
and a Toeplitz Fourier integral operator $\tilde{V}_{\chi}$ on
$\tilde{M}$ and a central operator $M_{\gamma}$ such that
$T_{\gamma}^* \tilde{V}_{\chi} T_{\gamma} = M_{\gamma}
\tilde{V}_{\chi}$, and such that $\alpha$ is induced by
conjugation by $\tilde{V}.$ But by definition, $\alpha$ is also
given by conjugation by  $U$. Now the Schwarz kernel of $U$, hence
$U$, lifts to a $\Gamma$ invariant kernel $\tilde{U}$ on
$\tilde{M}$. By assumption, $\tilde{U} A \tilde{U}^*$ has the same
complete symbol as $\tilde{V} A \tilde{V}^*$ for any Toeplitz
operator $A$ on $M$, lifted to ${\mathcal T}_{\Gamma}$.  By Lemma
\ref{E}, it follows that $\tilde{U} = \tilde{V} + R $, where $R$
is a smoothing Toeplitz operator. It follows that $M_{\gamma} = 1$
(hence $\theta_{\gamma, \chi} = 1$ for all $\gamma$),  and
therefore the symplectic map $\tilde{\chi}$ underlying
$\tilde{V}$, when lifted to $\tilde{X}$, is invariant under the
deck transformation group $\tilde{\Gamma}$ of $\tilde{X} \to X$.

\end{proof}

\section{\label{QTM} Quantization of torus maps}

To clarify the issues involved, we consider some standard examples
on the  symplectic 2m-torus ${\bf T}^{2m} = \C^m / \Z^{2nm}$.
Since $H_1({\bf T}^{2m}, \C) = \C^{2m}$, there will exist
symplectic maps which cannot be quantized in the sense of
Definitions \ref{TFIO} and \ref{TQM} (a) as quantum maps, though
they can and will be quantized as automorphisms. In fact, the
distinction can already be illustrated with the simplest maps:
\begin{itemize}

\item Kronecker translations $T_{\theta}(x) = x + \theta$ ($x, \theta \in {\bf T}^{2m}).$

\item Symplectic automorphisms $A \in Sp(2m, \Z).$

\end{itemize}

We begin by describing the line bundle on the torus and its
universal cover.  We follow \cite{Zel1, BSZ} and refer there for
further discussion.

The quotient setting is $L \to {\bf T}^{2n}$,  where $L$ is the
bundle with curvature $\sum_j dz_j \wedge d\bar{z}_j$. Sections of
$L^N$ are theta-functions of level $N$.  On the universal cover,
we have the pulled back bundle $L_\H = \C \times \C^m \to \C^m$.
Its associated principal $S^1$ bundle $\C^m \times S^1 \to \C^m$
is the {\it reduced Heisenberg group} ${\bf H}_{\rm red}^m$. We
recall that it is the quotient under the subgroup $(0, \Z)$ in the
center of the simply connected Heisenberg group ${\bf H}^m=\C^m
\times \R$ with  group law
$$(\zeta, t) \cdot (\eta, s) = (\zeta + \eta, t + s + \frac{1}{2} \Im (\zeta
\cdot \bar{\eta})).$$ The identity element is $(0, 0)$ and
$(\zeta, t)^{-1} = (- \zeta, - t)$. The  reduced Heisenberg group
is thus ${\bf H}^m_{\rm red}={\bf H}^m/ \{(0,  k): k \in \Z\} =
\C^m \times S^1$ with group law $$(\zeta, e^{2 \pi i t}) \cdot
(\eta, e^{2 \pi i s}) = (\zeta + \eta, e^{2 \pi i[t + s +
\frac{1}{2} \Im (\zeta \cdot \bar{\eta})]}).$$

We now equip ${\bf H}^m_{\rm red}$  with the  left-invariant
connection form \begin{equation} \alpha^L = \half  \sum_q(\xi_q \;
dx_q- x_q \;d\xi_q) - \frac{dt}{2 \pi},\;\; (\zeta=x+i \xi),
\end{equation} whose curvature equals  the symplectic form $\omega = \sum_q
dx_q \wedge d\xi_q$. The kernel of $\alpha^L$ is the distribution
of horizontal planes. To define the \szego kernel we further need
to split the complexified horizontal spaces into their holomorphic
and anti-holomorphic parts.  The left-invariant (CR-) holomorphic
(respectively anti-holomorphic) vector fields $Z_q^L$
(respectively $ \bar{Z}_q^L$) on ${\bf H}^m_{\rm red}$ are the
horizontal lifts  of the vector fields $\frac{\partial}{\partial
z_q},$
  respectively
$\frac{\partial}{\partial \bar z_q}$ with respect to $\alpha^L$.
We then define the Hardy space $\mathcal{H}^2 ({\bf H}^m_{\rm
red})$ of CR holomorphic functions to be the functions in
$\lcal^2({\bf H}^m_{\rm red})$ satisfying the left-invariant
Cauchy-Riemann equations $\bar{Z}^L_q f = 0$ ($1\le q\le m$) on
${\bf H}^m_{\rm red}$. For $N=1,2,\dots$, we further  define
$\hcal^2_N\subset \hcal^2(\H^m_{\rm red})$ as the
(infinite-dimensional) Hilbert space of square-integrable CR
functions $f$ such that $f\circ r_\theta = e^{iN\theta}f$ as
before. The representation $\mathcal{H}^2_1$ is irreducible and
may be identified with the Bargmann-Fock space of entire
holomorphic functions on $\C^n$ which are square integrable
relative to $e^{-|z|^2}$. The Szeg\"o kernel $\Pi_N^\H(x,y)$ is
the orthogonal projection to $\hcal^2_N$.  It is given by
\begin{equation}\label{szegoheisenberg} \Pi_N^\H(x,y)  =\frac{1}{\pi^m} N^m
e^{i N (t-s )} e^{ N(\zeta\cdot\bar \eta -\half |\zeta|^2
-\half|\eta|^2) } \,,\qquad x=(\zeta,t)\,,\
y=(\eta,s)\,.\end{equation} We note that it satisfies the
estimates in (\ref{OFFDIAG}). In this model example, Proposition
\ref{UPSZEGON}  was proven in \cite{Zel1}.

Finally, we describe the circle bundle $X$ in the quotient
setting. The lattice $\Z^{2m}$ may be embedded as a subgroup ${\bf
H}^m_{\Z}$ of ${\bf H}^m_{\R}$ under the homomorphism
\begin{equation} \iota(m,n) = (m, n, e^{i \pi  m \cdot n}). \end{equation}
We will denote the image by ${\bf H}^m_{\Z}$. To clarify the role
of the factors of $\half$ we show that ${\bf H}^m_{\Z}$ is indeed
a subgroup:
\begin{equation} \begin{array}{l} (m, n, e^{i \pi  m \cdot n}) \cdot (m', n', e^{i \pi  m' \cdot n'}) =
(m + m', n + n', e^{i \pi (m \cdot n + m' \cdot n' + m n' - m' n)})\\ \\
= (m + m', n + n', e^{i \pi (m \cdot n + m' \cdot n' + m n' + m' n)}) \\ \\
= (m + m', n + n', e^{i \pi ((m + m') \cdot (n + n'))}).
\end{array} \end{equation} We then put: $X = H^m_{\Z} \backslash H^nm_{\R}$,
i.e. $X$ is the {\it left} quotient of $H^m_{\R}$ by ${\bf
H}^m_{\Z}$. The left-invariant contact form $\alpha^L$ descends to
$X$ as a contact form and a connection form for the principal
$S^1$ bundle $X \to \C^m/\Z^{2m}$.

\subsection{Kronecker translations}

We first show that  irrational Kronecker translations
\begin{equation} T_{(a, b)} f( x, \xi) = f(x + a, \xi + b)\end{equation}
are non-quantizable as Toeplitz quantum maps.

\begin{prop}\label{KRON}  $T_{(a, b)}$ fails to be quantizable for all $(a, b) \in \R^{2n}/ \Z^{2n}.$
$T_{(a, b)}$ is quantizable at level $N$ iff $Na, Nb \in \Z^{2n}.$
\end{prop}

\begin{proof}

  By Proposition \ref{SYMPLIFT}, the
map lifts if and only if translations preserve holonomy of
homologically non-trivial loops. The loops on ${\bf T}^{2m}$ that
we need to consider are given in local coordinates by
$\gamma_{m,n}(t) = (t m, t n).$ Horizontal lifts to $X$ are given
by  $\tilde{\gamma}_{m,n}(t) = (tm, tn, 1)$. At $t = 1$ we obtain
$(m,n, 1) \sim (0,0, e^{i \pi m \cdot n}).$ Hence, the holonomy of
the path $\gamma_{m,n}$  equals $e^{i \pi m \cdot n}$.

Now translate the loop $\gamma_{(m,n)}$  by $(a, b)$ to obtain the
loop $\gamma_1(t) = (tm + a, tn + b)$. A  horizontal lift to $X$
is given by  $\tilde{\gamma}_1(t) = (tm + a, tn + b, e^{ \pi i (
b\cdot m - a \cdot n)t})$. It is the projection to $X$ of the left
translate by $(a, b, 0)$ of the original horizontal path. At $t =
1$ the endpoint is $(m + a, n + b, e^{ \pi i ( b\cdot m - a \cdot
n)}) = (a, b, e^{ i \pi [m \cdot n + 2 (b\cdot m - a \cdot n)] })
$. Hence, the holonomy changed by $e^{2 \pi i (b \cdot m - a \cdot
n)}$. The holonomy is preserved iff $b \cdot m - a \cdot n \in \Z$
for all $(m,n) \in \Z^{2m}$ iff $(a, b) \in \Z^{2m}$.

Lifting to level $N$ means changing the holonomy to $e^{2 \pi i N
\theta_{\gamma}}$. So the condition to lift becomes $(a, b ) \in
\frac{1}{N} \Z^2$

\end{proof}

\begin{rem} The  non-quantizability of $T_{(a, b)}$ is due to the left/right
invariance of various objects. $T_{(a,b)}$ only lifts to $H_{\Z}
\backslash H_{\R}$ as right translation by $(a, b, 0).$ But right
translation by an element of $H_{\R}$ does not preserve the left
invariant contact form. Equivalently, $T_{(a, b)}$ only lifts to a
contact transformation of $H_{\R}$ if it lifts to left translation
by $(a, b, 0)$. But then the lift does not descend to $H_{\Z}
\backslash H_{\R}$.
\end{rem}

\subsubsection{Kronecker translations as automorphisms}

It is easy to see that  Kronecker translations define
automorphisms of the revelant algebras.  We lift $T_{(a, b)}$ to
$H^n_{\R}$ as the contact transformation of left multiplication
$$T'_{(a, b)}(x, \xi, e^{2 \pi i t}) := (x + a, \xi + b, e^{2 \pi
i t} e^{ \pi i ( a \xi - b x)}).$$ Although the map $T_{(a, b)}$
does not descend to the quotient as a map, we claim:

\begin{prop} \label{KRON2} Kronecker maps $T_{(a, b)}$  have the following properties:

\begin{itemize}

\item (i) $T_{(a, b)}$  defines an automorphism of  ${\mathcal
T}^*_{\Gamma}$.

\item (ii)  $T_{(a, b)}$  defines an automorphism of ${\mathcal T}^{\infty}(X)/{\mathcal T}^{-\infty}(X)$ by
$$\alpha_{a, b; N}(Op_N(a) ) = Op_N(a \circ T_{(a, b)}).$$

\item (iii)  However, $\alpha_{(a, b)}$ does not preserve ${\mathcal K}_{\Gamma}$
and does not define an automorphism of ${\mathcal T}^{\infty}(X)$.

\end{itemize}
\medskip

\end{prop}

\begin{proof}

\noindent{\bf (i)} Left translation by $(a, b)$ defines  an
automorphism of ${\mathcal T}^*_{\Gamma}$ because, by
(\ref{szegoheisenberg}), the \szego kernel commutes with left
translations, i.e.
\begin{equation} \Pi_N^\H(\alpha \cdot x,\alpha \cdot y)  =
\Pi_N^\H( x, y)\;\;\; \forall \alpha. \end{equation} Indeed, it is
the kernel of a convolution operator.
\medskip

\noindent{\bf (ii)}

 Consider the conjugates $T_{\gamma}^{-1}
T_{a, b }' T_{\gamma}$ where $\gamma \in \Gamma = \Z^{2m}.$ An
easy computation shows that
$$T_{\gamma}^{-1} T_{a, b}' T_{\gamma} f(x, \xi, t) = M_{\gamma} T_{a, b}' f(x),\;\;\;\mbox{where}\;\;
 M_{\gamma}f(x, \xi, t) = f(x, \xi, t + \omega(\gamma,
(a, b)).$$ We need to show that $M_{\gamma}$ commutes with every
Toeplitz operator $\tilde{\Pi}_N \sigma \tilde{\Pi}_N$ with symbol
lifted from $\C^m/\Z^{2m}.$ Since the symbol is invariant under
the central circle, it is sufficient to show that $[M_{\gamma},
\tilde{\Pi}_N] = 0.$ But this follows as long as $\tilde{\Pi}_N (z
\cdot x, z \cdot y) = \tilde{\Pi}_N (x,y)$ for any $z = e^{i t}
\in S^1$, the center of the Heisenberg group. But this follows
because
$$\tilde{\Pi}_N (z \cdot x, z \cdot y) = |z|^{2N} \tilde{\Pi}_N (x,y) = \tilde{\Pi}_N (x,y).$$

Since left translation commutes with $\tilde{\Pi}_N$, the
automorphism  descends to the quotient as:

$$\begin{array}{lll} \alpha_N(\Pi_N a_N \Pi_N) & = &L^* T_{\alpha, \beta}^{' *} (\tilde{\Pi}_N
p^*a_N
 \tilde{\Pi}_N) T_{\alpha, \beta}'L \\ & & \\&
 = &  L^*\;  \tilde{\Pi}_N (T_{a, b}'\;  p^*a_N )\;  \tilde{\Pi}_N \;L. \end{array}$$
This is the stated formula.
\medskip

\noindent{\bf (iii)}  A  Kronecker automorphisms $\alpha_{(a, b)}$
can only preserve ${\mathcal K}_{\Gamma, N}$ if
 $\Pi_N a \Pi_N = 0$  implies
$\Pi_N (a\circ T_{a, b}) \Pi_N = 0$. But if this were the case,
the elements $\Pi_N e^{i \langle k, x \rangle} \Pi_N$ would  be
distinct  eigenoperators with eigenvalues $e^{i \langle k, (a, b)
\rangle}$. This contradicts the finite dimensionality of the
algebra for fixed $N$.

\end{proof}

\subsection{Quantum cat maps} We now show, in a similar way, that
symplectic linear maps of ${\bf T}^2$  always define quantum
automorphisms, even though they do not always define quantum maps.
We write  $g = \left( \begin{array}{ll} a & b \\ & \\ c & d
\end{array} \right) \in SL(2, \Z),$ and define $g (x, \xi) = (ax + b \xi, cx +
d \xi)$ on the torus. It lifts to the reduced Heisenberg group by
$g (x, \xi, t) = (g (x, \xi), t)$.

\begin{prop} \label{PARA} $T_g$
is quantizable  iff $a \cdot c, b \cdot d \in 2 \Z.$
\end{prop}

\begin{proof} We go through the same calculation as for Kronecker
translations. This time, the  horizontal lift of the transformed
loop is $(t(a \cdot m + b \cdot n), t (c \cdot m + d \cdot n),
1)$. At $t = 1$ the endpoint is
$$((a \cdot m + b \cdot n), (c\cdot  m + d\cdot n), 1) =  (0, 0, e^{i \pi (a \cdot m + b \cdot n)
\cdot (c \cdot m + d \cdot n)}).$$ Since the holonomy of the
original path was $e^{i \pi m \cdot n}$, the change in holonomy
equals
$$e^{i \pi (m \cdot n - (a \cdot m + b \cdot n)  (c \cdot m + d \cdot n))} = 1 \iff ac, bd  \in 2 \Z.$$
Here we use that $ad + bc \equiv 1 (\mbox{mod} 2 \Z). $
\end{proof}

\subsubsection{Linear maps as automorphisms}

It is known that quantizable linear maps (cat maps) on the
quotient define quantum maps with exact Egorov theorems (see e.g.
\cite{Zel1}). We now show that non-quantizable maps as well
defined automorphisms by the exact Egorov formula:

\begin{prop} \label{PARA2} $T_g$  defines  an automorphism of ${\mathcal T}^{\infty}(X)/{\mathcal T}^{-\infty}(X)$ by
$$\alpha_{g; N}(Op_N(a) ) = Op_N(a \circ T_g).$$\end{prop}

\begin{proof}
Consider the conjugates $T_{\gamma}^{-1} T_g T_{\gamma}$ where
$\gamma \in \Gamma.$ We have:
$$T_{\gamma}^{-1} T_g T_{\gamma} f(x, \xi, t) = M_{\gamma}T_g f(x),\;\;\;\mbox{where}\;\;
 M_{\gamma}f(x, \xi, t) = f((x, \xi) + (I - g) \gamma , t + \omega(\gamma, z)
 - \omega (\gamma, g(z + \gamma) ).$$
 $M_{\gamma}$ is the composition of translation $T_{(I -
 g)\gamma}$ with a central translation. Since $(I - g) \gamma $ is
 in the lattice, translation by this element commutes with left
 invariant operators. Thus, $M_{\gamma} \in {\mathcal
 T}_{\Gamma}^*(\tilde{X})'.$

Thus, the  automorphism  descends to the quotient as:

$$\begin{array}{lll} \alpha_N(\Pi_N a_N \Pi_N) & = &L^* T_g^* (\tilde{\Pi}_N
p^*a_N
 \tilde{\Pi}_N) T_gL \\ & & \\&
 = &  L^*\;  \tilde{\Pi}_N (T_g  p^*a_N )\;  \tilde{\Pi}_N \;L. \end{array}$$

\end{proof}

\section{\label{SP} Spectra of automorphisms}

In this article, our interest lies in the automorphisms defined by
symplectic maps. But  most  of the interest in quantizations of
quantizable  symplectic maps, at least in the physics literature,
is in their spectral theory as unitary operators $U_{\chi, N}$ on
the finite dimensional Hilbert spaces ${\mathcal H}_N (X)$. In
this section, we point out how the most important aspects of this
spectral theory of $U$ pertain only to the spectrum of the
associated  automorphism $U A U^*$. The main point is that the
reformulation suggests generalizations to other kinds of
automorphisms.  We also tie together the automorphisms of Toeplitz
algebras on the torus with the well-known ones on the rotation
algebra.

\subsection{Spectra of automorphisms of Hilbert-Schmidt algebras}
We  let ${\mathcal H}$ denote a  Hilbert space, and denote by
${\mathcal HS} $ the algebra of Hilbert-Schmidt  operators on
${\mathcal H}$, i.e. the operators for which the inner product
$\langle A, B \rangle:= Tr A B^*$ is finite.  We let $*$ denote
the adjoint  on ${\mathcal H}$. A finite dimensional algebra of
Hilbert-Schmidt operators is of course a full matrix algebra, and
its automorphisms are given by conjugation by unitary operators.

Suppose that $\alpha$ is an automorphism of ${\mathcal HS} $.

\begin{defn} We say that an automorphism $\alpha$ of $ {\mathcal HS}$  is:

\begin{itemize}

\item a * automorphism if: $\alpha(A^*) = \alpha(A)^*$.

\item {\em unitary} if: $\langle \alpha(A), \alpha(B) \rangle = \langle A, B \rangle.$

\item {\em a conjugation} if: there exists a unitary operator $ U: {\mathcal H} \to {\mathcal H}$ s.th. $\alpha(A) = U A U^*$.

\end{itemize}

\end{defn}

We also say that the automorphism is {\em tracial} if: $Tr
\alpha(A) = Tr A$ for all $A$ of trace class.

We will consider the eigenvalues and `eigenoperators' of a unitary
*-automorphism on ${\mathcal HS}$:
$$\alpha(A) = e^{i \theta} A.$$

If $\alpha$ is a unitary automorphism, then ( as a unitary
operator) it possesses an orthornormal basis of eigenoperators
$\{A_j\}$.

The following is elementary from the definitions.

\begin{prop} We have:

\noindent(i) A tracial *-automorphism is unitary.

\noindent(ii) The composition of any two eigenoperators is a
(possibly zero) eigenoperator.

\noindent(iii) If $A_j$ is an eigenoperator, then $A_j^* A_j$ is
an invariant operator, i.e. $\alpha(A_j^* A_j) = A_j^* A_j.$

\end{prop}

\begin{proof}

\noindent(i) Immediate from the fact that $$\langle \alpha(A),
\alpha(B) \rangle = Tr \alpha(A) \alpha(B)^*= Tr \alpha (AB^*) =
Tr AB^* = \langle A, B \rangle.$$
\medskip

\noindent(ii) - (iii) These  statements  follow from the equations
$\alpha(A_j A_k) = \alpha(A_j) \alpha(A_k) = e^{i (\theta_j +
\theta_k)}A_jA_k,$ and  $\alpha(A_j^*) = [\alpha(A_j)]^* = e^{- i
\theta_j} \alpha(A_j^*)$.

\end{proof}

In the case $\alpha (A) = U^* A U$ we note that the eigenoperators
of the automorphism  are given by \begin{equation}\label{EIGAUTO}
\alpha (\phi_{ N, j} \otimes \phi_{N, k}^*) = e^{i (\theta_{N ,j}
- \theta_{N, k})} \; (\phi_{N, j} \otimes \phi_{N,
k}^*),\end{equation} where $\{(\phi_{N, j}, e^{i \theta_{N,j}})\}$
are the spectral data of $U$.

\subsection{Spectral problems of quantum chaos}

The main problems on quantum maps  pertain to the spacings between
eigenvalues (the pair correlation problem) and the asymptotics of
matrix elements relative to eigenfunctions of the operators.

\subsubsection{Pair correlation problem}

Let us recall that the pair correlation function  $\rho_{2 N}$
(PCF) of a quantum map  $\{U_N\}$ with Planck constant $1/N$  is
the function  on $\R$ defined by
$$\int_{\R} f(x) d \rho_{2 N} (x) = \sum_{\ell = 0}^{\infty} \hat{f}(\frac{\ell}{N}) |Tr U_N^{\ell}|^2.$$
Its limit as $N \to \infty$, when one exists, is the PCF of the
quantum map.  Clearly, knowledge of $d \rho_2^N$ is equivalent to
knowledge of its {\it form factor} $|Tr U_N^{\ell}|^2$.  We
observe that the form factor depends only on the automorphism:

\begin{prop} \label{FORMFACTOR} The form factor of a conjugation automorphism $\{\alpha_N\}$ is given by
$$Tr \alpha_N^{\ell}  = \sum_j \langle \alpha_N^{\ell} A_j, A_j \rangle$$
where $\{A_j\}$ is an orthonormal basis for ${\mathcal HS}_N$.
\end{prop}

Indeed, if  $\alpha(A)_N = U_N^* A U_N$, then $\alpha_N^{\ell} =
U_N^{\ell} \otimes U_N^{-\ell}$ on ${\mathcal HS}_N$, and we have
$Tr_{{\mathcal HS}_N} U_N^{\ell} \otimes U_N^{-\ell} =
|Tr_{{\mathcal H}_N} U_N^{\ell} |^2.$

Hence,  the pair correlation problem makes sense for all unitary
automorphisms $\alpha_N$.
\medskip

\noindent{\bf Problem} Given any unitary  automorphism $\alpha_N$,
determine $\frac{1}{d_N} \sum_{j = 1}^{d_N^2} \delta(N
(\vartheta_{N,j}))$ as $N \to \infty,$ where
\begin{equation} \alpha_N(\Phi_j) = e^{i \vartheta_{N, j}} \Phi_j \end{equation}
is the eigenvalue problem for the automorphism.
\medskip

\subsubsection{Problems of quantum ergodicity/mixing}

We observe that these too can be formulated for any sequence of
automorphisms.  We rewrite the asymptotics of matrix elements
$\langle A \phi_{N, j}, \phi_{N, k} \rangle = Tr A \phi_{N, j}
\otimes \phi_{N, k}^*$ as the inner products $\langle A, \Phi_{jk}
\rangle$.  We observe that the eigenfunctions $\Phi_{k,k} =
\phi_{N, k} \otimes \phi_{N, k}^*$ always have eigenvalue $1$,
i.e. they are invariant states of the automorphism.

 It is simple to check that the proof of quantum
ergodicity for quantizations of ergodic quantizable  symplectic
maps $\chi$ (see \cite{Zel1}) uses only the automorphism involved.
It states that if a symplectic map $\chi$ is ergodic and
quantizable, then the invariant states of the corresponding
automorphism of the Toeplitz algebra are asymptotic to the traces
$\tau_N(A) = \frac{1}{\dim {\mathcal H}_M} Tr A |_{{\mathcal
H}_N}.$

It might be interesting to find generalizations of this result to
other kinds of automorphisms.

\subsection{Spectral theory of model automorphisms}

We now  point out that the automorphisms induced by model quantum
maps on the torus are  the same as the well-known automorphism of
the finite dimensional rotation algebras.

\subsubsection{The rotation algebra modulo $N$}

We denote by $G_{N} $ the finite Heisenberg group of order $N^2$,
generated by two  elements $U, V$ satisfying
$$U_2 U_1 = e^{2 \pi i/ N } U_1 U_2,\;\;\;U_1^N = U_2^N = I,$$
and its group algebra by  ${\mathcal R}_N$. $G_N$ has a unique
irreducible unitary  representation $\rho_N$ on $\C^N$ given as
follows: If we regard $\C^N = L^2(\Z / N \Z)$ then
$$\rho_N(U_1) \psi(Q) = e^{\frac{2 \pi i Q}{N}} \psi(Q),\;\;\; \rho_N(U_2) \psi(Q) = \psi(Q + 1).$$
 Recall that the  rotation algebra or   non-commutative torus $A_{\theta}$ is  the (pre-) C* algebra
generated by unitaries $U_1, U_2$ satisfying the Weyl commutation
relation
$$U_2 U_1 = e^{2 \pi i \alpha} U_1 U_2. $$
When $\alpha = 1/N$, $A_{\theta}$ has a large center generated by
$U_1^N, U_2^N$. ${\mathcal R}_N$ is obtained from $A_{\theta}$ by
viewing central elements as scalars.

\subsubsection{Toeplitz algebra and rotation algebra}

We  identify  the rotation algebra with $\theta = \frac{2\pi}{N}$
to the algebra of Toeplitz operators $\Pi_N a \Pi_N$ on the torus.
As verified by S. Nonnenmacher \cite{N}, the elements
\begin{equation} U_1 = e^{\pi N^2} \Pi_N e^{i \theta_1} \Pi_N,\;\;\; U_2 =
e^{\pi N^2} \Pi_N e^{i \theta_2} \Pi_N, \end{equation} satisfy
$U_j^N = I$.  Here, $(e^{i \theta_1}, e^{i \theta_2})$ are the
standard coordinates on the torus. Hence, any quantum map on the
torus  defines an automorphism of ${\mathcal R}_N$. Thus, we can
identify the automorphisms $\alpha_{g, N}$ quantizing  $g = \left(
\begin{array}{ll} a & b
\\ & \\ c & d
\end{array} \right)$ with the well-known automorphisms of
${\mathcal R}^N$ defined by $U_1 \to U_1^a U_2^b,\; U_2 \to U_1^c
U_2^d$ (see e.g. \cite{Narn}).

We can also see easily from this point of view that Kronecker maps
$T_{u,v}$  cannot in general  be quantized as automorphisms.
Namely, the quantization on $\R^2$ translates the symbol, so it
would descend to $U_1 \to e^{ i u} U_1, U_2 \to e^{iv} U_2$. To be
well-defined, one needs $e^{i u}, e^{i v}$ to be $N$th roots of
unity, which of course they are not in the irrational case.

\section{Appendix}

The key elements of the Toeplitz algebra and its automorphisms are
the $*$ product (\ref{STAR}) and the Egorov formula
(\ref{EGOROV}). The purpose of this appendix is to direct the
reader's attention to the existence of  routine calculations of
the complete symbols of compositions $a_N *_N b_N$ of symbols and
of conjugations $U_N Op_N(a) U_N^*$ of observables by Toeplitz
quantum maps. The method is to use the Boutet de Monvel -
Sj\"ostrand parametrix for the \szego kernel \cite{BS} as in
\cite{Zel2}.  Since the original version of this paper was
written, several papers \cite{KS, Sch1, Sch2, Sch3, STZ} have also
used this method to describe the $*$ product on symbols, so we
will be brief.

\begin{prop} Let $(M \omega)$ be a compact kahler manifold. Then:
 The  $*$ product defines an algebra structure on  classical symbols.  There exists an asymptotic expansion:
$$\hat{f}_1 * \hat{f}_2(z, \bar{z}) \sim \sum_{k = 0}^{\infty} N^{-k} B_k(f_1, f_2)$$
where $B_0(a,b) = f_1 \cdot f_2, B_1(f_2,f_2) = \half\{f_1,f_2\}$
and where $B_k$ is a bi-differential operator of $C^{\infty}(M)
\times C^{\infty}(M) \to C^{\infty}(M).$

\end{prop}

\begin{proof} Using the Boutet de Monvel-Sjostrand parametrix as
in \cite{Zel1, STZ}, one can obtain a complete asymptotic
expansion of the covariant symbol  $ \Pi_N a \Pi_N  b \Pi_N
(z,\bar{z}) $. One writes out $\Pi_N(z,w)$ as an oscillatory
integral and applies complex stationary phase. For calculations of
this kind we refer to \cite{STZ}.

To obtain $a *_N b$, we invert the Berezin transform $I_N$ on
symbols, as described in \cite{RT} and elsewhere. It is invertible
on formal power series, and the same inverse is well-defined on
symbol expansions. Thus,
$$a *_N b \sim  I_N^{-1} \Pi_N a \Pi_N  b \Pi_N
(z,\bar{z}).$$ This produces the symbol expansion claimed in the
proposition.

\end{proof}

\begin{prop} Let $U_{\chi, N}$ be a Toeplitz quantum map as in
Definition \ref{TQM}. Then for any observable $\Pi_N a_N \Pi_N =
Op_N(a_N)$, the contravariant symbol $a_{\chi}$ of $U_{\chi, N}^*
Op_N(a_N) U_{\chi, N} $ possesses a complete asymptotic expansion
$$a_{\chi} (z) \sim \sum_{k = 0}^{\infty} N^{-k} V_k (a \circ \chi) $$
where $V_0(a) = a$,  and where $V_k$ is a differential operator of
order at most $2k.$

\end{prop}

\begin{proof} As above, the expansion is obtained from the
covariant symbol $U_{\chi, N}^* Op_N(a_N) U_{\chi, N} (z,
\bar{z})$ by inverting the Berezin transform.  The asymptotics of
the covariant symbol follow by applying stationary phase to the
oscillatory integral formula for $U_{\chi, N}^* Op_N(a_N) U_{\chi,
N}$.

\end{proof}

\newpage

\end{document}